\documentclass[smallextended,referee,envcountsect,]{svjour3}
\smartqed
\usepackage{graphicx}
\usepackage{subfigure}

\journalname{}
\usepackage{mathptmx}      % use Times fonts if available on your TeX system
\usepackage{graphicx,epstopdf}
\usepackage{bm}
\usepackage{latexsym,enumerate}
\usepackage{mathrsfs}
\usepackage{amsfonts}
\usepackage{amssymb,amsmath}
\usepackage{epsfig}       %for Postscript files
\usepackage{subfigure}    %for parallel figures
\usepackage{graphicx}
\usepackage{float}
\usepackage{multirow}
\usepackage{array}
\usepackage{color}
\usepackage{listings}
\DeclareMathOperator*{\argmin}{argmin}
\usepackage[ruled]{algorithm2e}

\numberwithin{lemma}{section}
\numberwithin{theorem}{section}
\numberwithin{corollary}{section}
\numberwithin{definition}{section}
\numberwithin{property}{section}
\numberwithin{remark}{section}
\numberwithin{example}{section}
\numberwithin{figure}{section}
\numberwithin{proposition}{section}

\begin{document}

\title{Fast inertial dynamic algorithm with smoothing method for nonsmooth convex optimization}

\titlerunning{Fast inertial dynamic algorithm for nonsmooth convex optimization}
\author{Xin Qu         \and
        Wei Bian* %etc.
        }

\institute{\textbf{Statements and Declarations}
The authors declare that they have no conflict of interest.
\vspace{0.5cm}
             \\Xin Qu \at
             Harbin Institute of Technology \\
              School of Mathematics, Harbin 150001, P.R. China\\
              20B912001@stu.hit.edu.cn
           \and
              Wei Bian,  *Corresponding author  \at
              Harbin Institute of Technology \\
              School of Mathematics, Harbin 150001, P.R. China\\
              bianweilvse520@163.com
}

\date{Received: date / Accepted: date}
%The correct dates will be entered by the editor.

\maketitle

\begin{abstract}
In order to solve the minimization of a nonsmooth convex function, we design an inertial second-order dynamic algorithm, which is obtained by approximating the nonsmooth function by a class of smooth functions. By studying the asymptotic behavior of the dynamic algorithm, we prove that each trajectory of it weakly converges to an optimal solution under some appropriate conditions on the smoothing parameters, and the convergence rate of the objective function values is $o\left(t^{-2}\right)$. We also show that the algorithm is stable, that is, this dynamic algorithm with a perturbation term owns the same convergence properties when the perturbation term satisfies certain conditions. Finally, we verify the theoretical results by some numerical experiments.
\end{abstract}

%All acknowledgements should be placed in the back of the paper after Conclusions..

\section{Introduction}
Let $\mathcal{H}$ be a real Hilbert space endowed with the scalar product $\langle\cdot,\cdot\rangle$ and norm$ \|\cdot\|$. In this paper, our goal is to design an accelerated numerical method to solve the following convex optimization problem
\begin{equation}\label{x1}
\min_{x\in\mathcal{H}} f(x),
\end{equation}
where $f$ : $\mathcal{H}\rightarrow\mathbb{R}$ is a nonsmooth convex function. For the nonsmooth function $f$, we use a class of smooth convex functions to approximate it. Then we consider a dynamic algorithm with a smoothing function of $f$, that is,
\begin{equation}\label{prob2}
\ddot{x}(t)+\frac{\alpha}{t}\dot{x}(t)+\nabla_x\tilde{f}(x(t),\mu(t))= \textbf{0},
\end{equation}
where $\alpha>0$, $\tilde{f} : \mathcal{H} \times [0, \infty) \to \mathbb{R}$ is a smoothing function of convex function $f$, $\mu : [0,\infty)\rightarrow[0,\infty)$ is a continuously differentiable and decreasing function satisfying $\lim_{t\rightarrow\infty}\mu(t)=0$. The definition of smoothing function for convex function $f$ will be defined in Section 2. In our setting, because of the singularity of the damping coefficient $\frac{\alpha}{t}$ at $t=0$, we always set the initial time $t_0>0$.
Our main work is to study the asymptotic behavior of dynamic algorithm \eqref{prob2} for solving \eqref{x1}.
\subsection{Associated dynamic algorithms when $f$ is a smoothing function}
There is a long history of using dynamic algorithms to solve optimization problems \cite{ref45,ref46}. The asymptotic behavior of some dynamic algorithms has been studied when the function $f$ is smooth and convex.
The heavy ball with friction algorithm is one of them, which is modeled by
\begin{equation}\label{h1}
\ddot{x}(t)+\gamma\dot{x}(t)+\nabla f(x(t))=\textbf{0},
\end{equation}
where $\gamma$ is a fixed positive damping coefficient. This dynamic algorithm was first introduced by Polyak in \cite{ref10,ref9} from the perspective of optimization, and $\rm{\acute{A}}$lvarez studied the convergence of the  trajectories in the case of convexity in \cite{ref11}. For a general smooth convex function $f$, the convergence rate of dynamic algorithm \eqref{h1} is $O(t^{-1})$ in the worst case. When $f$ is strongly convex and $\gamma$ is selected appropriately, the convergence rate of dynamic algorithm \eqref{h1} can be exponential.
Since there is too much friction involved in this process, replacing the fixed viscosity coefficient with vanishing viscosity coefficient yields the inertial gradient dynamic algorithm
\begin{equation}\label{in2}
\ddot{x}(t)+\gamma(t)\dot{x}(t)+\nabla f(x(t))=\textbf{0},
\end{equation}
where $\gamma(\cdot)$ is a time-dependent positive damping coefficient. It has been studied by Cabot, Engler and Gaddat \cite{ref14,ref15}, and developed by Attouch and Cabot \cite{ref16}. A particularly interesting situation is the case $\gamma(t)=\frac{\alpha}{t}$. Su, Boyd and Cand$\rm{\grave{e}}$s \cite{Ref1} studied the following dynamic algorithm
\begin{equation}\label{in1}
\ddot{x}(t)+\frac{\alpha}{t}\dot{x}(t)+\nabla f(x(t))= \textbf{0},
\end{equation}
with $\alpha>0$. When $\alpha\geq3$, they proved that dynamic algorithm \eqref{in1} owns the fast convergence property $f(x(t))-\min f=O(t^{-2})$. Furthermore, \cite{ref2,ref5,ref3,ref4} showed that dynamic algorithm \eqref{in1} is the continuous version of the Nesterov accelerated gradient method with $\alpha=3$. Beck and Teboulle proposed the Fast Iterative Shrinkage Thresholding Algorithm (FISTA) in \cite{ref8} to solve the nonsmooth convex minimization problems with a splitting structure on the objective function, which is an extension of the accelerated gradient method in \cite{ref2}. Moreover, dynamic algorithm \eqref{in1} was further developed by Attouch-Chbani-Peypouquet-Redont to show that its each trajectory weakly converges to an element in $\argmin f$ when $\alpha>3$ \cite{ref6}. May \cite{ref7} proved that when $\alpha>3$, the asymptotic convergence rate of dynamic algorithm \eqref{in1} on the objective values can be improved from  $O(t^{-2})$ to $o(t^{-2})$. Chambolle-Dossal in \cite{ref17} had also obtained the same conclusions for the corresponding discrete algorithms. In the case of $\alpha\leq3$, Apidopoulos-Aujol-Dossal \cite{ref12} and Attouch-Chbani-Riahi \cite{ref13} demonstrated that the convergence rate of dynamic algorithm \eqref{in1} on the objective values is $O\left(t^{-\frac{2\alpha}{3}}\right)$. In addition, Attouch and Cabot \cite{ref20} studied the case that $f:\mathcal{H}\rightarrow\mathbb{R}\cup\{\infty\}$ is a convex lower semicontinuous proper function, and obtained the convergence rate on the objective values. The corresponding dynamic algorithm is
\begin{equation}\label{h2}
\ddot{x}(t)+\frac{\alpha}{t}\dot{x}(t)+\nabla f_{\lambda(t)}(x(t))=\textbf{0},
\end{equation}
where $f_\lambda:\mathcal{H}\rightarrow\mathbb{R}$ is the Moreau envelope of $f$ for index $\lambda>0$.
 Let us review some main properties of dynamic algorithm \eqref{h2}.
 \\$\bullet$  For $\alpha\geq3$, its trajectories satisfy the fast minimization property $f_{\lambda(t)}(x(t))-\min f=O\left(t^{-2}\right)$ and $f(\xi(t))-\min f=O\left(t^{-2}\right)$, where $\xi(t)={\rm prox}_{\lambda(t)f}(x(t))$ and ${\rm prox}_{\lambda f}(x)={\rm argmin}_{\zeta\in\mathcal{H}}\left\{f(\zeta)+\frac{1}{2\lambda}\|x-\zeta\|^2\right\}$.
 \\$\bullet$ For $\alpha>3$, the improved convergence rates are $f_{\lambda(t)}(x(t))-\min f=o\left(t^{-2}\right)$ and $f(\xi(t))-\min f=o\left(t^{-2}\right)$. In addition, each trajectory converges weakly to an optimal solution of $\min_{\mathcal{H}} f$ under appropriate conditions.
\subsection{Smoothing methods}
The subgradient methods were the first numerical schemes to solve nonsmooth convex minimization problems \cite{ref0}. It has been proved that the complexity of using these methods to obtain an $\varepsilon$-approximate solution of nonsmooth optimization problems is $O(\varepsilon^{-2})$. Smoothing methods are effective to overcome the nonsmoothness of optimization problems, which had been developed in the past decades \cite{ref27,refc,ref24,ref25}. Nesterov \cite{ref4} proposed a special smoothing technique for constructing efficient schemes for nonsmooth convex optimization, which is to approximate the initial nonsmooth objective function by a function with Lipschitz-continuous gradient. He showed that the complexity of finding an $\varepsilon$-approximate solution of nonsmooth optimization problems by smoothing technology is $O(\varepsilon^{-1})$. Chen introduced the smoothing methods for nonsmooth nonconvex minimization problems in \cite{ref26}, the main feature of which is to approximate nonsmooth functions by parameterized smoothing functions. She showed the properties of the smoothing functions and the gradient consistency of the subdifferentials related to a smoothing function, and presented how to update the smoothing parameter in the outer iteration of the smoothing methods to ensure that the iterative sequence converges to a stationary point of the original optimization problem.

These smoothing methods are widely used in various nonsmooth optimization problems. Zhang and Chen \cite{ref89} presented a novel smoothing active set method to solve the linearly constrained non-Lipschitz nonconvex minimization problems. They proved that any accumulation point of the iterative sequence generated by the smoothing active set method is a stationary point of the original problem. Bian and Chen studied the sparse regression problem with constraints in \cite{ref88}, the loss function of which is nonsmooth and convex. They gave an exact continuous relaxation model with the same optimal solution set as the regression problem, and then proposed a smoothing proximal gradient (SPG) algorithm based on the smoothing methods to find a lifted stationary point of the continuous relaxation problem. Burke-Chen-Sun \cite{bue} proposed an approximation theory of smooth functions for measurable composite max (CM) functions, explained the sub-consistency of gradient of CM integrands, and proved that the subgradient of expectation function can be approximated by smoothing without regularity.

On the one hand, we note that the dynamic algorithms \eqref{h1}-\eqref{in1} are not well-posed when $f$ is a convex lower semicontinuous proper function. On the other hand, though we can use the above second-order dynamic algorithm \eqref{h2} to solve the nonsmooth convex optimization problem \eqref{x1}, we need to know the Moreau envelope of $f$, which is much difficult for many functions. Thus, we will introduce smoothing methods into the dynamic algorithm, which is to use a sequence of smoothing  functions to approximate the nonsmooth function. The main advantage of the smoothing method is that we can easily construct the smoothing functions for a large class of nonsmooth functions. Thus, the dynamic algorithm not only can be well-defined, but also can be implemented easily.

This paper is organized as follows. In Section \ref{s2}, some preliminary results are presented, what's more, the existence and uniqueness of solutions of the considered dynamic algorithm \eqref{prob2} are proved. In Section \ref{s3}, we give the convergence rate on objective values $f(x(t))$ along the solution of dynamic algorithm \eqref{prob2}, and prove that the solution of it weakly converges to a minimizer of $f$. In Section \ref{s4}, we analyze the properties of dynamic algorithm \eqref{prob2} with a perturbation term. When the perturbation satisfies some appropriate conditions, the same convergence properties can be obtained. Finally, We use some numerical experiments to illustrate our theoretical results in Section \ref{s5}.
\section{Preliminaries}\label{s2}
For any $t\in\mathbb{R}$, we use $L^1(t,\infty)$ to denote the space of integrable functions from $(t,\infty)$ to $\mathbb{R}$, namely, $L^1(t,\infty):=\{\varphi:(t,\infty)\rightarrow\mathbb{R}|\int_{t}^{\infty}\varphi(s)ds<\infty\}$; $L^1_{loc}([t,\infty))$ denotes the space of locally integrable functions on $[t,\infty)$, that is, $L^1_{loc}([t,\infty)):=\{\varphi:[t,\infty)\rightarrow\mathbb{R}|\int_{t}^{T}\varphi(s)ds<\infty, \forall T\in{(t,\infty)}\}$. For a function $\phi:\mathcal{H}\rightarrow\mathbb{R}$, we let $[\phi]_+=\max\{\phi(x),0\}$, which is the positive part of function $\phi$.
\subsection{Smooth approximation}
A famous way to solve optimization problems with nonsmooth functions is to approximate these nonsmooth functions by a sequence of smooth functions.
This paper uses a class of smoothing functions defined as follows.
\begin{definition}\label{def1}\cite{ref88}
Let $f : \mathcal{H} \rightarrow \mathbb{R}$ be a convex function. We call $\tilde{f} : \mathcal{H} \times [0, \infty) \to \mathbb{R}$ a smoothing function of $f$, if $\tilde{f}(x, \mu)$ satisfies the following conditions:
\begin{enumerate}[{\rm (i)}]
\item for any fixed $\mu > 0$, $\tilde{f}(\cdot, \mu)$ is continuously differentiable in $\mathcal{H}$, and for any fixed $x \in \mathcal{H}$, $\tilde{f}(x, \cdot)$ is continuously differentiable in $(0, \infty)$;

\item $\lim_{z\to x, \mu \downarrow 0}\tilde{f}(z,\mu)=f(x), \quad \forall x\in \mathcal{H}$;

\item $\tilde{f}(x,\mu)$ is convex with respect to $x$ in $\mathcal{H}$ for any fixed $\mu>0$;

\item there exists a positive constant $\kappa$ such that
\begin{equation}\label{fka}
|\nabla_\mu\tilde{f}(x,\mu)|\leq \kappa, \quad \forall\mu\in(0,\infty), x\in \mathcal{H};
\end{equation}
\item there exists a constant $L>0$ such that for any $\mu \in(0, \infty)$, $\nabla_x \tilde{f}(\cdot,\mu)$ is Lipschitz continuous on $\mathcal{H}$ with Lipschitz constant $L\mu^{-1}$;
\item $\nabla_z \tilde{f}(z,\mu)$ is continuous with respect to $\mu$ on $(0,\infty)$ for any fixed $z\in\mathcal{H}$.
\end{enumerate}
\end{definition}
By Definition \ref{def1}-(iv), we know that
\begin{equation}\label{fk}
|\tilde{f}(x,\mu_2)-\tilde{f}(x,\mu_1)|\leq \kappa|\mu_1-\mu_2|, \quad \forall x\in \mathcal{H},\ \mu_1, \mu_2 \in (0, \infty).
\end{equation}
Furthermore,
\begin{equation}\label{fkb}
|\tilde{f}(x,\mu)-f(x)|\leq \kappa \mu, \quad \forall x\in \mathcal{H},\  \mu\in(0,\infty).
\end{equation}
For the function $\mu(\cdot)$ in dynamic algorithm \eqref{prob2}, the following hypothesis is assumed throughout the paper:
\begin{equation*}
 \bm{(H_1)} : \int_{t_0}^{\infty}t\mu(t)dt<\infty.
\end{equation*}
\begin{remark}
For fixed $t_0>0$, $(H_1)$ implies
\begin{equation}\label{ass}
\int_{t_0}^{\infty} \frac{1}{t}\mu(t)dt<\infty.
\end{equation}
\end{remark}
\subsection{Preliminary results}
Before giving the existence and uniqueness of the global solution to \eqref{prob2}, we introduce some lemmas which are used in the following to analyze the asymptotic behavior of trajectories.
\begin{lemma}\label{ap1}{\rm\cite{ref28}}
 Let $S$ be a nonempty subset of $\mathcal{H}$ and let $x : [0,\infty)\rightarrow\mathcal{H}$. Assume that
\begin{enumerate}[{\rm (i)}]
\item for every $z\in S$, $\lim_{t\rightarrow\infty}\|x(t)-z\|$ exists;
\item every weak sequential limit point of $x(t)$, as $t\rightarrow\infty$, belongs to $S$.
\end{enumerate}
Then $x(t)$ converges weakly as $t\rightarrow\infty$ to a point in $S$.
\end{lemma}
\begin{lemma}\label{ap2}{\rm\cite{ref29}}
Take $\delta>0$, and let $g\in L^1(\delta,\infty)$ be nonnegative and continuous. Consider a nondecreasing function $\psi : [\delta,\infty)\rightarrow[0,\infty)$ such that $\lim_{t\rightarrow\infty}\psi(t)=\infty$. Then
\begin{equation*}
\lim_{t\rightarrow\infty}\frac{1}{\psi(t)}\int_{\delta}^{t}\psi(s)g(s)ds=0.
\end{equation*}
\end{lemma}
\begin{lemma}\label{ap3}{\rm{\cite{ref6}}}
 Let $\delta>0$, and $w : [\delta,\infty)\rightarrow\mathbb{R}$ be a continuously differentiable function which is bounded from below. Assume
\begin{equation*}
t\ddot{w}(t)+\alpha\dot{w}(t)\leq m(t),
\end{equation*}
for some $\alpha>1$, almost every $t>\delta$, and some nonnegative function $m\in L^1(\delta,\infty)$.
Then, $[\dot{w}]_+\in L^1(t_0,\infty)$, and $\lim_{t\rightarrow\infty}w(t)$ exists.
\end{lemma}
\begin{lemma}\label{ap4}{\rm\cite{ref30}}\underline{}
Let $m : [\delta,T]\rightarrow[0,\infty)$ be integrable, and constant $c\geq0$. Suppose $w : [\delta,T]\rightarrow\mathbb{R}$ is continuous and
\begin{equation*}
\frac{1}{2}w^2(t)\leq\frac{1}{2}c^2+\int_{\delta}^{t}m(\tau)w(\tau)d\tau,
\end{equation*}
for all $t\in[\delta,T]$. Then, $|w(t)|\leq c+\int_{\delta}^{t}m(\tau)d\tau$ for all $t \in [\delta,T]$.
\end{lemma}
\subsection{Existence and uniqueness of solutions}
\begin{proposition}\label{propo1}
For every initial value $x_0:=x(t_0)\in\mathcal{H}$ and $v_0:=\dot{x}(t_0)\in\mathcal{H}$, there exists a unique global trajectory
$x : [t_0, \infty)\rightarrow\mathcal{H}$ of the dynamic algorithm \eqref{prob2}.
\end{proposition}
{\it Proof}\quad Denote $X(t):=\left(
\begin{array}{c}
x(t)\\
\dot{x}(t)
\end{array}
\right)$ and let $F:[ t_0,\infty)\times\mathcal{H}\times\mathcal{H}\rightarrow \mathcal{H}\times\mathcal{H}$ be
\[F(t,z,v)=
\left(
\begin{array}{c}
v\\
-\frac{\alpha}{t}v-\nabla_z\tilde{f}(z,\mu(t))
\end{array}
\right).
\]
We endow $\mathcal{H}\times\mathcal{H}$ with scalar product $\left\langle(z,v),(\bar{z},\bar{v})\right\rangle_{\mathcal{H}\times\mathcal{H}}=\langle z,\bar{z}\rangle+\langle v,\bar{v}\rangle$ and norm $\|(z,v)\|_{\mathcal{H}\times\mathcal{H}}=\|z\|+\|v\|$. Hence \eqref{prob2} can be written as
\begin{eqnarray}\label{xgs18}
\left\{\aligned
  & \frac{d}{dt}X(t)=F(t, X(t)),
  \\&  X(t_0)=\left(
  \begin{array}{c}
  x_0\\
  v_0
  \end{array}
  \right).
\endaligned\right.
\end{eqnarray}
For the first-order dynamic algorithm \eqref{xgs18}, we apply the non-autonomous version of Cauchy-Lipschitz-Picard theorem \cite{ref33} to prove the existence and uniqueness of its solution.
\\ \textbf{Step1:} For any $(z,v)$, $(\bar{z},\bar{v})\in\mathcal{H}\times\mathcal{H}$, from Definition \ref{def1}-(v), we know that there exists $L>0$ such that
\begin{equation*}
\begin{split}
& \|F(t,z,v)-F(t,\bar{z},\bar{v})\|_{\mathcal{H}\times\mathcal{H}}
\\= &  \|v-\bar{v}\|+\left\|-\frac{\alpha}{t}(v-\bar{v})+\nabla_z\tilde{f}(\bar{z},\mu(t))-\nabla_z\tilde{f}(z,\mu(t))\right\|
\\ \leq & \left(1+\frac{\alpha}{t}\right)\|v-\bar{v}\|+L\mu(t)^{-1}\|z-\bar{z}\|
\\ \leq & M(t)\left\|(z,v)-(\bar{z}-\bar{v})\right\|_{\mathcal{H}\times\mathcal{H}},
\end{split}
\end{equation*}
where $M(t)=\max\left\{\left(1+\frac{\alpha}{t}\right), L\mu(t)^{-1}\right\}$,$\quad\forall t\in[t_0, \infty)$. Hence $F(t,\cdot,\cdot)$ is $M(t)$-Lipschitz continuous for every $t\geq t_0$. Moreover, for any $t\geq t_0$, by the continuity of $\frac{\alpha}{t}$ and $\mu(t)^{-1}$, we know that $M(\cdot)$ is integrable on $[t_0,T]$ for any $t_0<T<\infty$. Thus $M(\cdot)\in L_{loc}^1\left([t_0,\infty)\right)$.
\\ \textbf{Step2:} For fixed $z,v\in\mathcal{H}$, $t_0<T<\infty$, we get
\begin{equation*}
\begin{split}
& \int_{t_0}^{T}\|F(t,z,v)\|_{\mathcal{H}\times\mathcal{H}}dt
\\=& \int_{t_0}^{T}\left(\|v\|+\left\|-\frac{\alpha}{t}v-\nabla_z\tilde{f}(z,\mu(t))\right\|\right)dt
\\\leq & \int_{t_0}^{T}\left(\left(1+\frac{\alpha}{t}\right)\|v\|+\left\|\nabla_z\tilde{f}(z,\mu(t))\right\|\right)dt.
\end{split}
\end{equation*}
By Definition \ref{def1}-(vi) and the continuity of $\mu(\cdot)$, we know that $\nabla_z\tilde{f}(z,\mu(t))$ is continuous with respect to $t$ for fixed $z$. This together with the continuity of $\frac{\alpha}{t}$ yields
$$\int_{t_0}^{T}\|F(t,z,v)\|_{\mathcal{H}\times\mathcal{H}}dt<\infty,\quad\forall t_0<T<\infty.$$
\textbf{Step3:} For fixed $z,v\in\mathcal{H}$, we obtain
\begin{equation}\label{mm2}
\begin{split}
\|F(t,z,v)\|_{\mathcal{H}\times\mathcal{H}} & =\|v\|+\left\|-\frac{\alpha}{t}v-\nabla_z\tilde{f}(z,\mu(t))\right\|
\\& \leq\left(1+\frac{\alpha}{t}\right)\|v\|+\left\|\nabla_z\tilde{f}(z,\mu(t))\right\|.
\end{split}
\end{equation}
In view of Definition \ref{def1}-(v), we know
\begin{equation}\label{mm1}
\left\|\nabla_z\tilde{f}(z,\mu(t))\right\|\leq\left\|\nabla_z\tilde{f}(\textbf{0},\mu(t))\right\|+L\mu(t)^{-1}\|z\|,\quad\forall z\in\mathcal{H}.
\end{equation}
Substituting \eqref{mm1} into \eqref{mm2}, we get
\begin{equation}
\begin{split}
 \|F(t,z,v)\|_{\mathcal{H}\times\mathcal{H}}& \leq\left(1+\frac{\alpha}{t}\right)\|v\|+\left\|\nabla_z\tilde{f}(\textbf{0},\mu(t))\right\|+L\mu(t)^{-1}\|z\|
\\& \leq P(t)\left(1+\|v\|+\|z\|\right),
\end{split}
\end{equation}
where $P(t):=1+\frac{\alpha}{t}+L\mu(t)^{-1}+\left\|\nabla_z\tilde{f}(\textbf{0},\mu(t))\right\|$. By virtue of the continuity of $\frac{\alpha}{t}$, $\mu(t)$ and $\nabla_z\tilde{f}(\textbf{0},\mu(t))$ with respect to $t$, we conclude that $P(t)\in L_{loc}^1(t_0,\infty)$.
Therefore, by  Cauchy-Lipschitz-Picard theorem, we can obtain that there is a global unique solution for dynamic algorithm \eqref{xgs18}, and then the proof is completed.
\qed
\section{Convergence of dynamic algorithm in \eqref{prob2}}\label{s3}
In this section, we will analyze the convergence properties of trajectory to dynamic algorithm \eqref{prob2}, including the convergence rate on the objective values and the weak convergence of the trajectory to a minimizer of $f$.
\subsection{Minimizing property}
We begin by introducing a function that plays a crucial role in proving weak convergence of the trajectory to \eqref{prob2}.
\\Let $z\in\mathcal{H}$, and define the function $h: (t_0,\infty)\rightarrow\mathbb{R^+}$ by
\begin{equation}\label{num1}
h(t)=\frac{1}{2}\|x(t)-z\|^2.
\end{equation}
By differentiating it, we obtain\\
$$\dot{h}(t)=\langle x(t)-z, \dot{x}(t)\rangle  \quad\mbox{and}\quad \ddot{h}(t)=\langle x(t)-z, \ddot{x}(t)\rangle+\|\dot{x}(t)\|^2.$$
Using \eqref{prob2} and the convex inequality of $\tilde{f}(x,\mu)$ with respect to $x$, we have
\begin{equation*}
\begin{split}
\ddot{h}(t)+\frac{\alpha}{t}\dot{h}(t)= & \|\dot{x}(t)\|^2+\langle x(t)-z,-\nabla_x\tilde{f}(x(t),\mu(t))\rangle
\\ \leq & \|\dot{x}(t)\|^2+\tilde{f}(z,\mu(t))-\tilde{f}(x(t),\mu(t)).
\end{split}
\end{equation*}
Rearranging the terms, we find
\begin{equation}\label{q5}
\ddot{h}(t)+\frac{\alpha}{t}\dot{h}(t)+\tilde{f}(x(t),\mu(t))-\tilde{f}(z,\mu(t))\leq\|\dot{x}(t)\|^2.
\end{equation}
\begin{proposition}\label{soul}
Suppose $\alpha>0$ and $\inf f>-\infty$. Let $x : [t_0,\infty)\rightarrow\mathcal{H}$ be a trajectory of \eqref{prob2}. Then $$ \lim_{t\rightarrow\infty}f(x(t))=\inf f,$$
and $$\lim_{t\rightarrow\infty}\|\dot{x}(t)\|=0.$$
\end{proposition}
{\it Proof}\quad
In view of \eqref{fkb}, we know
\begin{equation}\label{in4}
\tilde{f}(z,\mu(t))\leq\kappa\mu(t)+f(z),\quad \forall z\in \mathcal{H}.
\end{equation}
 Introducing \eqref{in4} into \eqref{q5}, we conclude that
\begin{equation}\label{q6}
\ddot{h}(t)+\frac{\alpha}{t}\dot{h}(t)+\tilde{f}(x(t),\mu(t))-f(z)-\kappa\mu(t)\leq\|\dot{x}(t)\|^2.
\end{equation}
Let us introduce the function $W : [t_0,\infty)\rightarrow\mathbb{R}$ defined by
$$W(t):=\frac{1}{2}\|\dot{x}(t)\|^2+\tilde{f}(x(t),\mu(t))+\kappa\mu(t).$$
Differentiating $W$ along the trajectory of \eqref{prob2} and by \eqref{fka}, we obtain
\begin{equation}\label{q7}
\frac{d}{dt}W(t)=-\frac{\alpha}{t}\|\dot{x}(t)\|^2+\left(\nabla_\mu\tilde{f}(x(t),\mu(t))+\kappa\right)\dot{\mu}(t)\leq-\frac{\alpha}{t}\|\dot{x}(t)\|^2\leq0.
\end{equation}
Thus the function $W$ is nonincreasing on $[t_0,\infty)$. Recalling $\inf f>-\infty$ and \eqref{fkb}, we have
 \begin{equation*}
 W(t)\geq\tilde{f}(x(t),\mu(t))+\kappa\mu(t)\geq f(x(t))\geq\inf f>-\infty.
 \end{equation*}
Hence $W_\infty=\lim_{t\rightarrow\infty}$W(t) exists.
From \eqref{q7}, it follows that
\begin{equation}\label{q8}
\int_{t_0}^{\infty}\frac{1}{t}\|\dot{x}(t)\|^2dt\leq\frac{1}{\alpha}\left(W(t_0)-W_\infty\right)<\infty.
\end{equation}
\\Substituting $W$ in \eqref{q6}, we get
\begin{equation}\label{q9}
\ddot{h}(t)+\frac{\alpha}{t}\dot{h}(t)+W(t)-f(z)\leq b(t),
\end{equation}
where $b(t):=\frac{3}{2}\|\dot{x}(t)\|^2+2\kappa\mu(t)$.
Multiplying each member of inequality \eqref{q9} by $t^\alpha$, we find
\begin{equation*}
\frac{d}{dt}\left(t^\alpha\dot{h}(t)\right)+ t^\alpha\left(W(t)-f(z)\right)\leq t^\alpha b(t).
\end{equation*}
Integrating the above inequality on $[t_0,t]$, we obtain
\begin{equation*}
t^\alpha\dot{h}(t)-{t^\alpha_0}\dot{h}(t_0)+\int_{t_0}^{t}s^\alpha(W(s)-f(z))ds\leq\int_{t_0}^{t}s^\alpha b(s)ds.
\end{equation*}
By virtue of the nonincreasing property of $W$, we deduce that
\begin{equation*}
t^\alpha\dot{h}(t)-{t^\alpha_0}\dot{h}(t_0)+(W(t)-f(z))\int_{t_0}^{t}s^\alpha ds\leq\int_{t_0}^{t}s^\alpha b(s) ds.
\end{equation*}
Dividing the above inequality by $t^\alpha$ and integrating it from $t_0$ to $t$ again, we get
\begin{equation*}
\begin{split}
& h(t)-h(t_0)+\int_{t_0}^{t}s^{-\alpha}\left(W(s)-f(z)\right)\left(\int_{t_0}^{s}\tau^\alpha d\tau\right)ds
\\ \leq & t^\alpha_0\dot{h}(t_0)\int_{t_0}^{t}s^{-\alpha}ds+\int_{t_0}^{t}s^{-\alpha}\left(\int_{t_0}^{s}\tau^\alpha b(\tau)d\tau\right)ds.
\end{split}
\end{equation*}
Since $W$ is nonincreasing, we find
\begin{equation*}
\begin{split}
& h(t)-h(t_0)+(W(t)-f(z))\int_{t_0}^{t}s^{-\alpha}\left(\int_{t_0}^{s}\tau^\alpha d\tau\right)ds
\\ \leq & \frac{1}{\alpha-1}t_0\left|\dot{h}(t_0)\right|
+\int_{t_0}^{t}s^{-\alpha}\left(\int_{t_0}^{s}\tau^\alpha b(\tau)d\tau\right)ds.
\end{split}
\end{equation*}
Calculating and rearranging the above inequalities, we get
\begin{equation}\label{q10}
\begin{split}
& \frac{1}{\alpha+1}(W(t)-f(z))\left(\frac{t^2}{2}-\frac{t^2_0}{2}+\frac{t^{\alpha+1}_0}{(\alpha-1)t^{\alpha-1}}-\frac{t^2_0}{\alpha-1}\right)
\\\leq & h(t_0)+\frac{1}{\alpha-1}t_0\left|\dot{h}(t_0)\right|+\int_{t_0}^{t}s^{-\alpha}\left(\int_{t_0}^{s}\tau^\alpha b(\tau)d\tau\right)ds.
\end{split}
\end{equation}
By using Fubini theorem to estimate the last term in \eqref{q10}, we have
\begin{equation*}
\begin{split}
\int_{t_0}^{t}s^{-\alpha}\left(\int_{t_0}^{s}\tau^\alpha b(\tau)d\tau\right)ds= & \int_{t_0}^{t}\left(\int_{\tau}^{t}s^{-\alpha}ds\right)\tau^\alpha b(\tau)d\tau
\\ = & \frac{1}{\alpha-1}\int_{t_0}^{t}\left(\frac{1}{\tau^{\alpha-1}}-\frac{1}{t^{\alpha-1}}\right)\tau^\alpha b(\tau)d\tau
\\ \leq & \frac{1}{\alpha-1}\int_{t_0}^{t}\tau b(\tau)d\tau.
\end{split}
\end{equation*}
Coming back to inequality \eqref{q10}, we conclude that
\begin{equation*}
\begin{split}
& \frac{1}{\alpha+1}\left(W(t)-f(z)\right)\left(\frac{t^2}{2}-\frac{t^2_0}{2}+\frac{t^{\alpha+1}_0}{(\alpha-1)t^{\alpha-1}}-\frac{t^2_0}{\alpha-1}\right)
\\ \leq & h(t_0)+\frac{1}{\alpha-1}t_0\left|\dot{h}(t_0)\right|+\frac{1}{\alpha-1}\int_{t_0}^{t}\tau b(\tau)d\tau.
\end{split}
\end{equation*}
Dividing the above inequality by $t^2$ and rewriting the last term, we deduce that
\begin{equation}\label{q11}
\begin{split}
& \frac{1}{\alpha+1}(W(t)-f(z))\left(\frac{1}{2}-\frac{t^2_0}{2t^2}+\frac{t^{\alpha+1}_0}{(\alpha-1)t^{\alpha+1}}-\frac{t^2_0}{(\alpha-1)t^2}\right)
\\ \leq & \frac{h(t_0)+\frac{1}{\alpha-1}t_0\left|\dot{h}(t_0)\right|}{t^2}+\frac{1}{(\alpha-1)t^2}\int_{t_0}^{t}\tau^2 \frac{1}{\tau} b(\tau)d\tau.
\end{split}
\end{equation}
Under condition $\eqref{ass}$ and estimation \eqref{q8}, we have
\begin{equation*}
\int_{t_0}^{\infty}\frac{1}{t}b(t)dt=\int_{t_0}^{\infty}\left(\frac{3}{2t}\|\dot{x}(t)\|^2+2\kappa\frac{1}{t}\mu(t)\right)dt<\infty.
\end{equation*}
Taking the limit as $t\rightarrow\infty$ in \eqref{q11} and applying Lemma \ref{ap2}, we derive that
\begin{equation}\label{q12}
\limsup_{t\rightarrow\infty}W(t)\leq f(z).
\end{equation}
Under the definition of $W$ and by \eqref{fkb}, we know
\begin{equation*}
\limsup_{t\rightarrow\infty}f(x(t))\leq\limsup_{t\rightarrow\infty}\left(\tilde{f}(x(t),\mu(t))+\kappa\mu(t)\right)\leq f(z).
\end{equation*}
Since the above inequality holds for an arbitrary $z$, we conclude that
\begin{equation*}
\limsup_{t\rightarrow\infty}f(x(t))\leq \inf f.
\end{equation*}
Thus, we have
\begin{equation*}
\lim_{t\rightarrow\infty}f(x(t))=\inf f.
\end{equation*}
Then, \eqref{q12} and $W(t)\geq\frac{1}{2}\|\dot{x}(t)\|^2+f(x(t))$ further implies that
\begin{equation*}
\lim_{t\rightarrow\infty}\|\dot{x}(t)\|=0.
\end{equation*}
\qed
\subsection{Convergence rate on objective values}
\begin{theorem}\label{th}
Let $x : [t_0,\infty)\rightarrow\mathcal{H}$ be a trajectory of \eqref{prob2}, and assume {\rm argmin}$f$ is nonempty and $(H_1)$ is ture.
\begin{enumerate}[{\rm (i)}]
\item Suppose $\alpha\geq 3$. Then
\begin{equation*}
f(x(t))-{\rm min} f={O}\left(\frac{1} {t^2}\right).
\end{equation*}
\item Suppose $\alpha> 3$. Then
\begin{equation}\label{fkc}
\int_{t_0}^{\infty}t\left(f(x(t)-{\rm min} f\right)dt<\infty,\\
\end{equation}
\begin{equation}\label{fkd}
\int_{t_0}^{\infty}t\|\dot{x}(t)\|^2dt<\infty,
\end{equation}
\begin{equation}\label{fke}
f(x(t))-{\rm min} f= o\left(\frac{1} {t^2}\right).
\end{equation}
\end{enumerate}
\end{theorem}
{\it Proof}\quad
(i) Fix $x^* \in \textrm{argmin} f$, and consider the energy function
\begin{equation*}
\mathcal{E}(t)=t^{2}\left(\tilde{f}(x(t),\mu(t))-\tilde{f}(x^*,\mu(t))+2\kappa\mu(t)\right)+\frac{1}{2}\|(\alpha-1)\left(x(t)-x^*\right)+t\dot{x}(t)\|^2.
\end{equation*}
In view of \eqref{fkb},  this gives
\begin{equation}\label{fkf}
\tilde{f}(x(t),\mu(t))-\tilde{f}(x^*,\mu(t))+2\kappa\mu(t)\geq f(x(t))-f(x^*)\geq0,
\end{equation}
which implies $\mathcal{E}(t)\geq0$.\\
Using the classical derivation chain rule and equation \eqref{prob2}, we obtain
\begin{equation}\label{in3}
\begin{split}
\frac{d}{dt}\mathcal{E}(t)= & 2t\left(\tilde{f}(x(t),\mu(t))-\tilde{f}(x^*,\mu(t))+2\kappa\mu(t)\right) \\ & +t^2\left(\left\langle\nabla_x\tilde{f}(x(t),\mu(t)),\dot{x}(t)\right\rangle+\nabla_{\mu}\tilde{f}(x(t),\mu(t))\dot{\mu}(t)-\nabla_{\mu}\tilde{f}(x^*,\mu(t))\dot{\mu}(t)\right.
\\ & \left.+2\kappa\dot{\mu}(t)\right)+\left\langle(\alpha-1)(x(t)-x^*)+t\dot{x}(t),\alpha\dot{x}(t)+t\ddot{x}(t)\right\rangle
\\= & 2t\left(\tilde{f}(x(t),\mu(t))-\tilde{f}(x^*,\mu(t))+2\kappa\mu(t)\right)
\\ & -(\alpha-1)t\langle x(t)-x^*,\nabla_x\tilde{f}(x(t),\mu(t))\rangle
\\ & +t^2\left(\nabla_\mu\tilde{f}(x(t),\mu(t))\dot{\mu}(t)-\nabla_\mu\tilde{f}(x^*,\mu(t))\dot{\mu}(t)+2\kappa \dot{\mu}(t)\right).
\end{split}
\end{equation}
By \eqref{fka} and $\dot{\mu}(t)\leq0$, we deduce that
\begin{equation}\label{fkg}
\nabla_\mu\tilde{f}(x(t),\mu(t))\dot{\mu}(t)-\nabla_\mu\tilde{f}(x^*,\mu(t))\dot{\mu}(t)\leq -2\kappa\dot{\mu}(t).
\end{equation}\label{}
Since $\tilde{f}(x,\mu)$ is convex with respect to $x$ for any fixed $\mu$ , we have\\
\begin{equation}\label{fkh}
\tilde{f}(x^*,\mu(t))-\tilde{f}\left(x(t),\mu(t)\right)\geq\left\langle\nabla_x\tilde{f}(x(t),\mu(t)),x^*-x(t)\right\rangle.
\end{equation}
\\If $\alpha\geq 3$, introducing \eqref{fkg} and \eqref{fkh} into \eqref{in3}, we obtain
\begin{equation}\label{fkj}
\begin{split}
\frac{d}{dt}\mathcal{E}(t)\leq & -(\alpha-3)t\left(\tilde{f}(x(t),\mu(t))-\tilde{f}(x^*,\mu(t))+2\kappa\mu(t)\right)+2(\alpha-1)\kappa t\mu(t)\\
                  \leq & 2(\alpha-1)\kappa t\mu(t),
\end{split}
\end{equation}
where last inequality uses \eqref{fkf}.
 \\From $(H_1)$, we have the positive part $[\frac{d}{dt}\mathcal{E}(t)]_+$ belongs to $L_1(t_0,\infty)$. Let $\rho(t):=\mathcal{E}(t)-\int_{t_0}^{t}[\frac{d}{dt}\mathcal{E}(s)]_+ds$, and $\rho(\cdot)$ is bounded by the boundedness of $\mathcal{E}(\cdot)$ and $[\frac{d}{dt}\mathcal{E}(t)]_+
 \\\in L_1(t_0,\infty)$. This together with $\frac{d}{dt}\rho(t)=\frac{d}{dt}\mathcal{E}(t)-[\frac{d}{dt}\mathcal{E}(t)]_+\leq0$ yields the existence of $\lim_{t\rightarrow\infty}\rho(t)$. Hence,
 \begin{equation}\label{op}
 \lim_{t\rightarrow\infty}\mathcal{E}(t)=\lim_{t\rightarrow\infty}w(t)+\int_{t_0}^{\infty}[\mathcal{E}(s)]_+ds<\infty.
 \end{equation}
 It ensues that $\mathcal{E}(\cdot)$ is bounded on $[t_0,\infty)$, and then $t^2 \left(\tilde{f}(x(t),\mu(t))-\tilde{f}(x^*,\mu(t))+2\kappa\right.
 \\\left.\mu(t)\right)$ is bounded on $[t_0,\infty)$, which means that there exists $C>0$ such that
 $$t^2 \left(\tilde{f}(x(t),\mu(t))-\tilde{f}(x^*,\mu(t))+2\kappa\mu(t)\right)\leq C,\quad\forall t\geq t_0.$$
  As a consequence, returning to \eqref{fkf}, we have
 \begin{equation*}
 f(x(t))-f(x^*)\leq \frac{C}{t^2},
 \end{equation*}
namely
  \begin{equation*}
 f(x(t))-f(x^*)={O}\left(\frac{1} {t^2}\right).
 \end{equation*}
(ii) Now suppose $\alpha>3$. By integrating \eqref{fkj} from $t_0$ to $t$, we obtain
\begin{equation*}
\begin{split}
& \int_{t_0}^{t}s\left(\tilde{f}(x(s),\mu(s))-\tilde{f}(x^*,\mu(s))+2\kappa\mu(s)\right)ds
\\ \leq & \frac{1}{\alpha-3}\left(\mathcal{E}(t_0)-\mathcal{E}(t)\right)+\frac{2(\alpha-1)\kappa}{\alpha-3}\int_{t_0}^{t}s\mu(s)ds\\
\leq & \frac{1}{\alpha-3}\mathcal{E}(t_0)+\frac{2(\alpha-1)\kappa}{\alpha-3}\int_{t_0}^{t}s\mu(s)ds.
\end{split}
\end{equation*}
Under $(H_1)$, we have the estimate
\begin{equation}\label{q1}
\int_{t_0}^{\infty}s\left(\tilde{f}(x(s),\mu(s))-\tilde{f}(x^*,\mu(s))+2\kappa\mu(s)\right)ds<\infty.
\end{equation}
From \eqref{fkf}, we obtain
\begin{equation*}
\int_{t_0}^{\infty}t\left(f(x(t)-{\rm min} f\right)dt<\infty,
\end{equation*}
which shows \eqref{fkc}.
\\To prove \eqref{fkd}, take the scalar product of \eqref{prob2} with $t^2\dot{x}(t)$, then we have
\begin{equation*}
\begin{split}
& \frac{t^2}{2}\frac{d}{dt}\|\dot{x}(t)\|^2+\alpha t \|\dot{x}(t)\|^2+t^2\frac{d}{dt}\left(\tilde{f}(x(t),\mu(t))-\tilde{f}(x^*,\mu(t))+2\kappa\mu(t)\right)
\\= & t^2\left(\nabla_\mu\tilde{f}(x(t),\mu(t))\dot{\mu}(t)-\nabla_\mu\tilde{f}(x^*,\mu(t))\dot{\mu}(t)+2\kappa\dot{\mu}(t)\right).
\end{split}
\end{equation*}
By integrating the above equation from $t_0$ to $t$, we obtain
\begin{equation*}
\begin{split}
& \frac{t^2}{2}\|\dot{x}(t)\|^2-\frac{{t_0}^2}{2}\|\dot{x}(t_0)\|^2+(\alpha-1)\int_{t_0}^{t}s\|\dot{x}(s)\|^2ds+t^2\left(\tilde{f}(x(t),\mu(t))\right.
\\ & \left.-\tilde{f}(x^*,\mu(t))+2\kappa\mu(t)\right)-t_0^2\left(\right.\tilde{f}(x(t_0),\mu(t_0))-\tilde{f}(x^*,\mu(t_0))
\\ & +2\kappa\mu(t_0))-2\int_{t_0}^{t}s\left(\tilde{f}(x(s),\mu(s))-\tilde{f}(x^*,\mu(s))+2\kappa\mu(s)\right)ds
\\= & \int_{t_0}^{t}s^2\left(\nabla_\mu\tilde{f}(x(s),\mu(s))\dot{\mu}(s)-\nabla_\mu\tilde{f}(x^*,\mu(s))\dot{\mu}(s)+2\kappa\dot{\mu}(s)\right)ds.
\end{split}
\end{equation*}
Combining the above relation with \eqref{fkf}, \eqref{fkg} and $\frac{t^2}{2}\|\dot{x}(t)\|^2\geq0$, we conclude that
\begin{equation*}
\begin{split}
(\alpha-1)\int_{t_0}^{t}s\|\dot{x}(s)\|^2ds\leq\tilde{C}+2\int_{t_0}^{t}s\left(\tilde{f}(x(s),\mu(s))-\tilde{f}(x^*,\mu(s))+2\kappa\mu(s)\right)ds.
\end{split}
\end{equation*}
where $\tilde{C}=\frac{{t_0}^2}{2}\|\dot{x}(t_0)\|^2+t_0^2\left(\tilde{f}(x(t_0),\mu(t_0))-\tilde{f}(x^*,\mu(t_0))+2\kappa\mu(t_0)\right).$ By virtue of \eqref{q1} and $\alpha>3$, we get \eqref{fkd}.
\\Now, we consider
\begin{equation}\label{num}
E(t):=\frac{1}{2}\|\dot{x}(t)\|^2+\tilde{f}(x(t),\mu(t))-\tilde{f}(x^*,\mu(t))+2\kappa\mu(t),
\end{equation}
which is nonnegative on $[t_0,\infty)$. By \eqref{fkd} and \eqref{q1}, we know
\begin{equation}\label{q4}
\int_{t_0}^{\infty}tE(t)dt<\infty.
\end{equation}
\\Differentiating $t^2E$, we get that
\begin{equation*}
\begin{split}
\frac{d}{dt}(t^2E(t))= & 2tE(t)+t^2\frac{dE(t)}{dt}
\\= & 2t\left(\frac{1}{2}\|\dot{x}(t)\|^2+\tilde{f}(x(t),\mu(t))-\tilde{f}(x^*,\mu(t))+2\kappa\mu(t)\right)
\\ & +t^2(\langle\dot{x}(t),\ddot{x}(t)\rangle+\langle\nabla_x\tilde{f}(x(t),\mu(t)),\dot{x}(t)\rangle
\\& +\nabla_\mu\tilde{f}(x(t),\mu(t))\dot{\mu}(t) -\nabla_\mu\tilde{f}(x^*,\mu(t))\dot{\mu}(t)+2\kappa\dot{\mu}(t)).
\end{split}
\end{equation*}
By \eqref{prob2}, \eqref{fkg} and $\alpha>3$, we deduce that
\begin{equation}\label{in8}
\begin{split}
\frac{d}{dt}(t^2E(t))\leq & 2t\left(\tilde{f}(x(t),\mu(t))-\tilde{f}(x^*,\mu(t))+2\kappa \mu(t)\right)+t(1-\alpha)\|\dot{x}(t)\|^2
\\ \leq & 2t\left(\tilde{f}(x(t),\mu(t))-\tilde{f}(x^*,\mu(t))+2\kappa \mu(t)\right).
\end{split}
\end{equation}
Combining \eqref{q1} and \eqref{in8} yields that $[\frac{d}{dt}(t^2E(t))]_+\in$ $L_1(t_0,\infty)$, hence similar to the analysis for \eqref{op}, we know $\lim_{t\rightarrow\infty}t^2E(t)$ exists. Recalling \eqref{q4}, we have $\int_{t_0}^{\infty}tE(t)dt=\int_{t_0}^{\infty}\frac{1}{t}(t^2E(t))dt<\infty$. By $\int_{t_0}^{\infty}\frac{1}{t}dt=\infty$ and the existence of  $\lim_{t\rightarrow\infty}t^2E(t)$, we conclude that
\begin{equation*}
\lim_{t\rightarrow\infty}t^2E(t)=0.
\end{equation*}
Under the definition of $E$, we have the estimate
$$0\leq\lim_{t\rightarrow\infty}t^2\left(\tilde{f}(x(t),\mu(t))-\tilde{f}(x^*,\mu(t))+2\kappa\mu(t)\right)\leq\lim_{t\rightarrow\infty}t^2E(t)=0.$$
Furthermore,
$$0\leq \lim_{t\rightarrow\infty}t^2\left(f(x(t))-f(x^*)\right)\leq\lim_{t\rightarrow\infty}t^2\left(\tilde{f}(x(t),\mu(t))-\tilde{f}(x^*,\mu(t))+2\kappa\mu(t)\right)=0,$$
namely,
$$f(x(t))-{\rm min} f= o\left(\frac{1} {t^2}\right).$$
\qed
\subsection{Weak convergence of trajectories}
\begin{theorem}\label{theo}
Suppose {\rm argmin}$f \neq\emptyset$ and $(H_1)$ is ture. Let $x : [t_0,\infty)\rightarrow\mathcal{H}$ be the trajectory of \eqref{prob2} with $\alpha>3$, then $x(t)$ converges weakly in $\mathcal{H}$, as $t\rightarrow\infty$, to a point in {\rm argmin} $f$.
\end{theorem}
{\it Proof}\quad
The proof is based on the Opial's lemma (Lemma $\ref{ap1}$). For any $ x^* \in {\rm argmin} f,$ coming back to \eqref{q5}, and let $z=x^*$, we have
\begin{equation*}
t\ddot{h}(t)+\alpha\dot{h}(t)+t\left(\tilde{f}(x(t),\mu(t))-\tilde{f}(x^*,\mu(t))+2\kappa\mu(t)\right)\leq k(t),
\end{equation*}
where $k(t):=t\|\dot{x}(t)\|^2+2\kappa t\mu(t)$. From \eqref{fkf}, we have
\begin{equation*}
t\ddot{h}(t)+\alpha\dot{h}(t)\leq k(t).
\end{equation*}
Combining \eqref{fkd} with $(H_1)$, we know $k\in L^1(t_0,\infty)$. By applying Lemma \ref{ap3}, we know $\lim_{t\rightarrow\infty}h(t)$ exists. The first point of Opial's lemma is proved. It also implies that the trajectory $x(\cdot)$ is bounded on $[t_0,\infty)$. The next step is to prove the second point of Opail's lemma, which is that every weak sequential limit point of $x(t)$ belongs to {\rm argmin}$f$. Let $\bar{x}$ be a sequential limit point of $x(\cdot)$ on $[t_0,\infty)$ with convergence sequence $\{t_n\}$. In view of Proposition \ref{soul}, we have
\begin{equation*}
f(\bar{x})=\lim_{n\rightarrow\infty}f(x(t_n))=\lim_{t\rightarrow\infty}f(x(t))=\inf f.
\end{equation*}
It implies $\bar{x}\in{\rm argmin}f$, which gives the claim.
\qed
\section{Asymptotic convergence of \eqref{prob2} for minimization under perturbations}\label{s4}
In this section, we show that the perturbation term satisfying certain conditions does not affect the convergence results of \eqref{prob2} for solving optimization problem \eqref{x1}, that is, dynamic algorithm $\eqref{prob2}$ is stable. For this purpose, we consider the following dynamic algorithm with perturbation
\begin{equation}\label{prob3}
\ddot{x}(t)+\frac{\alpha}{t}\dot{x}(t)+\nabla_x\tilde{f}(x(t),\mu(t))=g(t),
\end{equation}
where $\alpha>0$, $g:[t_0,\infty)\rightarrow\mathcal{H}$ is the perturbation term and the functions $\tilde{f},\mu$ are defined same as in \eqref{prob2}.
\\For the function $g(\cdot)$ in dynamic algorithm \eqref{prob3}, the following hypothesis is assumed
throughout the paper:
\begin{equation*}
\bm{(H_g)}: \int_{t_0}^{\infty}\|g(s)\|ds<\infty.
\end{equation*}
Under the condition $(H_g)$ and Definition \ref{def1}-(v)(vi), we can prove the global existence and uniqueness of trajectory to \eqref{prob3} by similar analogy with Proposition \ref{propo1}.
\subsection{Minimizing property under perturbations}
\begin{proposition}\label{soul1}
Suppose $\alpha>0$, $\inf f>-\infty$ and $(H_g)$ is ture. Let $x : [t_0,\infty)\rightarrow\mathcal{H}$ be a trajectory of \eqref{prob3}. Then $$ \lim_{t\rightarrow\infty}f(x(t))=\inf f.$$
\end{proposition}
{\it Proof}\quad
Let $T>t_0$, and $t_0\leq t\leq T$. Define the energy function $W_{g,T}: [t_0,\infty)\rightarrow\mathbb{R}$ by
\begin{equation*}
W_{g,T}(t):=\frac{1}{2}\|\dot{x}(t)\|^2+\tilde{f}(x(t),\mu(t))-\inf f+\kappa\mu(t)+\int_{t}^{T}\langle\dot{x}(\tau),g(\tau)\rangle d\tau.
\end{equation*}
Differentiating $W_{g,T}$ and using \eqref{fka}, \eqref{prob3}, we get that
\begin{equation*}
\begin{split}
\frac{d}{dt}W_{g,T}(t)= & \langle\dot{x}(t),\ddot{x}(t)+\nabla_x\tilde{f}(x(t),\mu(t))-g(t)\rangle+\left(\nabla_\mu\tilde{f}(x(t),\mu(t))+\kappa\right)\dot{\mu}(t)\\ \leq & -\frac{\alpha}{t}\|\dot{x}(t)\|^2.
\end{split}
\end{equation*}
Hence $W_{g,T}$ is a nonincreasing function on $[t_0,\infty)$, which means that $W_{g,T}(t)\leq W_{g,T}(t_0)$, for any $t\geq t_0$, i.e.
\begin{equation*}
\begin{split}
&\frac{1}{2}\|\dot{x}(t)\|^2+\tilde{f}(x(t),\mu(t))-\inf f+\kappa\mu(t)+\int_{t}^{T}\langle\dot{x}(\tau),g(\tau)\rangle d\tau
\\\leq & \frac{1}{2}\|\dot{x}(t_0)\|^2+\tilde{f}(x(t_0),\mu(t_0))-\inf f+\kappa\mu(t_0)+\int_{t_0}^{T}\langle\dot{x}(\tau),g(\tau)\rangle d\tau.
\end{split}
\end{equation*}
From \eqref{fkb}, we know that $\tilde{f}(x(t),\mu(t))-\inf f+\kappa\mu(t)\geq0$.  Combining with Cauchy-Schwarz inequality, we obtain
\begin{equation*}
\frac{1}{2}\|\dot{x}(t)\|^2\leq\frac{1}{2}\|\dot{x}(t_0)\|^2+\tilde{f}(x(t_0),\mu(t_0))-\inf f+\kappa\mu(t_0)+\int_{t_0}^{t}\|\dot{x}(\tau)\|\|g(\tau)\| d\tau.
\end{equation*}
Together Lemma \ref{ap4} with $(H_g)$, we deduce that
\begin{equation*}
\|\dot{x}(t)\|\leq\left(\|\dot{x}(t_0)\|^2+2\left(\tilde{f}(x(t_0),\mu(t_0))-\inf f+\kappa\mu(t_0)\right)\right)^{\frac{1}{2}}+\int_{t_0}^{t}\|g(\tau)\| d\tau,
\end{equation*}
which gives
\begin{equation}\label{q14}
\sup_{t\geq t_0}\|\dot{x}(t)\|<\infty.
\end{equation}
From \eqref{q14}, we know that
\begin{equation*}
W_g(t):=\frac{1}{2}\|\dot{x}(t)\|^2+\tilde{f}(x(t),\mu(t))-\inf f+\kappa\mu(t)+\int_{t}^{\infty}\langle\dot{x}(\tau),g(\tau)\rangle d\tau
\end{equation*}
is well defined, and
\begin{equation*}
W_g(t)\geq-\sup_{t\geq0}\|\dot{x}(t)\|\int_{t_0}^{\infty}\|g(\tau)\|d\tau,
\end{equation*}
\begin{equation}\label{q15}
\frac{d}{dt}W_g(t)= \frac{d}{dt}W_{g,T}(t)\leq-\frac{\alpha}{t}\|\dot{x}(t)\|^2,
\end{equation}
which gives $\lim_{t\rightarrow\infty}W_g(t)=W_{g,\infty}\in\mathbb{R}.$
 \\Let us now integrate inequality \eqref{q15} on $[t_0,t]$ and let $t\rightarrow\infty$, we find
 \begin{equation}\label{q16}
 \int_{t_0}^{\infty}\frac{\alpha}{\tau}\|\dot{x}(\tau)\|^2d\tau\leq W_g(t_0)-W_{g,\infty}<\infty.
 \end{equation}
 For any $z\in\mathcal{H}$, recalling the function $h(t)=\frac{1}{2}\|x(t)-z\|^2$, similar to the analysis for \eqref{q5} and using \eqref{prob3}, we obtain
 \begin{equation}\label{q17}
 \ddot{h}(t)+\frac{\alpha}{t}\dot{h}(t)+\tilde{f}(x(t),\mu(t))-\tilde{f}(z,\mu(t))\leq\|\dot{x}(t)\|^2+\langle g(t),x(t)-z\rangle.
 \end{equation}
 By virtue of \eqref{fkb}, we get
 \begin{equation}\label{q18}
  \ddot{h}(t)+\frac{\alpha}{t}\dot{h}(t)+\tilde{f}(x(t),\mu(t))-f(z)-\kappa\mu(t)\leq\|\dot{x}(t)\|^2+\langle g(t),x(t)-z\rangle.
 \end{equation}
Substituting $W_g$ into \eqref{q18} and using Cauchy-Schwarz inequality to get
\begin{equation*}
\begin{split}
W_{g,\infty}+\inf f-f(z)\leq & W_{g}(t)+\inf f-f(z)\leq \frac{3}{2} \|\dot{x}(t)\|^2+\|g(t)\|\|x(t)-z\|
  \\+ & \sup_{t\geq t_0}\|\dot{x}(t)\|\int_{t}^{\infty}\|g(\tau)\|d\tau+2\kappa\mu(t)-\frac{1}{t^\alpha}\frac{d}{dt}\left(t^\alpha\dot{h}(t)\right),
\end{split}
\end{equation*}
where the first inequality uses the nonincreasing property of $W_g$.
\\Let $\theta>t_0.$ Multiplying both sides of the above inequality by $\frac{1}{t}$ and integrating it from $t_0$ to $\theta$, we conclude that
\begin{equation}\label{q19}
\begin{split}
\left(W_{g,\infty}+\inf f-f(z)\right)\ln\frac{\theta}{t_0}\leq & \frac{3}{2}\int_{t_0}^{\theta}\frac{1}{t}\|\dot{x}(t)\|^2dt+\int_{t_0}^{\theta}\frac{\|g(t)\|\|x(t)-z\|}{t}dt
\\ & + \sup_{t\geq t_0}\|\dot{x}(t)\|\int_{t_0}^{\theta}\left(\frac{1}{t}\int_{t}^{\infty}\|g(\tau)\|d\tau\right)dt+2\kappa\int_{t_0}^{\theta}\frac{1}{t}\mu(t)dt
\\ & -\int_{t_0}^{\theta}\frac{1}{t^{\alpha+1}}
\frac{d}{dt}\left(t^\alpha\dot{h}(t)\right)dt.
\end{split}
\end{equation}
Estimation of several terms in \eqref{q19} are given below.
\begin{enumerate}[{\rm (1)}]
\item Since
\begin{equation}\label{qwe}
\|x(t)-z\|\leq\|x(t_0)-z\|+\int_{t_0}^{t}\|\dot{x}(s)\|ds,
\end{equation}
we get
\begin{equation*}
\int_{t_0}^{\theta}\frac{\|g(t)\|\|x(t)-z\|}{t}dt\leq\left(\frac{\|x(t_0)-z\|}{t_0}+\sup_{t\geq t_0}\|\dot{x}(t)\|\right)\int_{t_0}^{\theta}\|g(t)\|dt,
\end{equation*}
where we use $\frac{1}{t}\int_{t_0}^{t}\|\dot{x}(s)\|ds\leq\sup_{t\geq t_0}\|\dot{x}(t)\|.$
\item By direct calculating, we have
\begin{equation*}
\begin{split}
\int_{t_0}^{\theta}\left(\frac{1}{t}\int_{t}^{\infty}\|g(\tau)\|d\tau\right)dt= & \int_{t_0}^{\theta}\left(\int_{t}^{\infty}\|g(\tau)\|d\tau\right)d\left(\ln t\right)
\\=& \ln\theta\int_{\theta}^{\infty}\|g(\tau)\|d\tau-\ln t_0\int_{t_0}^{\infty}\|g(\tau)\|d\tau
\\ & +\int_{t_0}^{\theta}\|g(t)\|\ln tdt.
\end{split}
\end{equation*}
\item Estimating the last term in inequality \eqref{q19}, we deduce that
\begin{equation*}
\begin{split}
\int_{t_0}^{\theta}\frac{1}{t^{\alpha+1}}\frac{d}{dt}\left(t^\alpha\dot{h}(t)\right)dt= & \frac{1}{\theta}\dot{h}(\theta)-\frac{1}{t_0}\dot{h}(t_0)-\int_{t_0}^{\theta}t^\alpha\dot{h}(t)(-\alpha-1)t^{-\alpha-2}dt
\\= & \frac{1}{\theta}\dot{h}(\theta)-\frac{1}{t_0}\dot{h}(t_0)
\\ & +(\alpha+1)\left(\frac{1}{\theta^2}h(\theta)-\frac{1}{t_0^2}h(t_0)+2\int_{t_0}^{\theta}\frac{1}{t^3}h(t)dt\right)
 \\= & C_0+\frac{1}{\theta}\dot{h}(\theta)+\frac{\alpha+1}{\theta^2}h(\theta)+2(\alpha+1)\int_{t_0}^{\theta}\frac{1}{t^3}h(t)dt
\\\geq & C_0+\frac{1}{\theta}\dot{h}(\theta),
\end{split}
\end{equation*}
where $C_0=(\alpha+1)\frac{1}{t_0^2}h(t_0)-\frac{1}{t_0}\dot{h}(t_0)$, and by \eqref{qwe}, we have
\begin{equation*}
|\dot{h}(\theta)| =\left|\langle\dot{x}(\theta), x(\theta)-z\rangle\right|\leq\sup_{t\geq t_0}\|\dot{x}(t)\|\left(\|x(t_0)-z\|+\theta\sup_{t\geq t_0}\|\dot{x}(t)\|\right).
\end{equation*}
\end{enumerate}
Combining the above results with \eqref{ass} and \eqref{q19}, we obtain
\begin{equation*}
\begin{split}
\left(W_{g,\infty}+\inf f-f(z)\right)\ln\frac{\theta}{t_0}\leq & C_1+\sup_{t\geq t_0}\|\dot{x}(t)\|\ln\theta\int_{\theta}^{\infty}\|g(t)\|dt
\\ & + \left(\sup_{t\geq t_0}\|\dot{x}(t)\|\right)\int_{t_0}^{\theta}\|g(t)\|\ln tdt,
\end{split}
\end{equation*}
where $C_1=\frac{3}{2}\int_{t_0}^{\theta}\frac{1}{t}\|\dot{x}(t)\|^2dt+\left(\frac{\|x(t_0)-z\|}{t_0}+\sup_{t\geq t_0}\|\dot{x}(t)\|\right)\int_{t_0}^{\theta}\|g(t)\|dt-\int_{t_0}^{\infty}\|g(t)\|dt
\\\ln t_0\sup_{t\geq t_0}\|\dot{x}(t)\|+2\kappa\int_{t_0}^{\theta}\frac{1}{t}\mu(t)dt+\frac{1}{\theta}\sup_{t\geq t_0}\|\dot{x}(t)\|\left(\|x(t_0)-z\|+\theta\sup_{t\geq t_0}\|\dot{x}(t)\|\right)$ is a constant by (1), (2) and (3). Dividing both sides by $\ln\frac{\theta}{t_0}$, letting $\theta\rightarrow\infty$ in the above inequality, and using Lemma \ref{ap2}, we have $W_{g,\infty}\leq f(z)-\inf f$,$\quad \forall z\in \mathcal{H}$, which implies $W_{g,\infty}\leq0$ by the continuity and convexity of $f$.
\\In fact
\begin{equation*}
 W_g(t)\geq\tilde{f}(x(t),\mu(t))-\inf f+\kappa\mu(t)-\sup_{t\geq t_0}\|\dot{x}(t)\|\int_{t}^{\infty}\|g(\tau)\|d\tau.
\end{equation*}
Letting $t\rightarrow\infty$, we obtain
\begin{equation*}
0\geq W_{g,\infty}\geq\limsup_{t\rightarrow\infty}\left(\tilde{f}(x(t),\mu(t))-\inf f+\kappa\mu(t)\right)\geq0,
\end{equation*}
which implies
\begin{equation*}
\lim_{t\rightarrow\infty}f(x(t))=\lim_{t\rightarrow\infty}\left(\tilde{f}(x(t),\mu(t))+\kappa\mu(t)\right)=\inf f.
\end{equation*}
\qed
\subsection{Convergence rate on the objective values under perturbations}
As can be seen from the following theorem, the convergence rate of the objective values along the trajectory of \eqref{prob3} is consistent with \eqref{prob2} under the condition of $\int_{t_0}^{\infty}t\|g(t)\|
\\dt<\infty$.
\begin{remark}
It is clear that $\int_{t_0}^{\infty}t\|g(t)\|dt<\infty$ implies $(H_g)$.
\end{remark}
\begin{theorem}\label{h5}
Let ${\rm argmin} f\neq\emptyset$ and $x:[t_0,\infty)\rightarrow\mathcal{H}$ be the trajectory of \eqref{prob3}. Assume $\int_{t_0}^{\infty}t\|g(t)\|dt<\infty$ and $(H_1)$ is ture.
\begin{enumerate}[{\rm (i)}]
\item If $\alpha\geq 3$. Then
\begin{equation*}
f(x(t))-{\rm min} f={O}\left(\frac{1} {t^2}\right).
\end{equation*}
\item If $\alpha> 3$. Then $x$ is bounded on $[t_0,\infty)$, and
\begin{equation}\label{q30}
\int_{t_0}^{\infty}t\left(f(x(t)-{\rm min} f\right)dt<\infty,
\end{equation}
\begin{equation}\label{q20}
\int_{t_0}^{\infty}t\|\dot{x}(t)\|^2dt<\infty,
\end{equation}
\begin{equation}\label{q21}
f(x(t))-{\rm min} f= o\left(\frac{1} {t^2}\right).
\end{equation}
\end{enumerate}
\end{theorem}
{\it Proof}\quad
(i) Given $x^*\in \argmin f$ and $T>t_0$,  we introduce the energy function $t\mapsto\mathcal{E}_{g,T}$ defined by
\begin{equation*}
\begin{split}
\mathcal{E}_{g,T}(t):= & \frac{2}{\alpha-1}t^2\left(\tilde{f}(x(t),\mu(t))-\tilde{f}(x^*,\mu(t))+2\kappa\mu(t)\right)
\\ & +(\alpha-1)\left\|x(t)-x^*+\frac{t}{\alpha-1}\dot{x}(t)\right\|^2\\ & +2\int_{t}^{T}\tau\left\langle x(\tau)-x^*+\frac{\tau}{\alpha-1}\dot{x}(t),g(\tau)\right\rangle d\tau.
\end{split}
\end{equation*}
By differentiating $\mathcal{E}_{g,T}$, and together \eqref{fkg} with \eqref{prob3}, we immediately find
\begin{equation}\label{ap5}
\begin{split}
\frac{d}{dt}\mathcal{E}_{g,T}(t)= & \frac{4}{\alpha-1}t(\tilde{f}(x(t),\mu(t))-\tilde{f}(x^*,\mu(t))+2\kappa\mu(t))
\\ & +\frac{2}{\alpha-1}t^2(\left\langle\nabla_x\tilde{f}(x(t), \mu(t)),\dot{x}(t)\right\rangle+\nabla_\mu\tilde{f}(x(t),\mu(t))\dot{\mu}(t)
\\ & - \nabla_\mu\tilde{f}(x^*,\mu(t))\dot{\mu}(t)+2\kappa\dot{\mu}(t))
\\ & +2(\alpha-1)\left\langle x(t)-x^*+\frac{t}{\alpha-1}\dot{x}(t),\frac{t}{\alpha-1}\left(\frac{\alpha}{t}\dot{x}(t)+\ddot{x}(t)-g(t)\right)\right\rangle
\\ \leq & \frac{4}{\alpha-1}t\left(\tilde{f}(x(t),\mu(t))-\tilde{f}(x^*,\mu(t))+2\kappa\mu(t)\right)
\\ & +\frac{2}{\alpha-1}t^2\left\langle\nabla_x\tilde{f}(x(t),\mu(t)),\dot{x}(t)\right\rangle
\\ & +2(\alpha-1)\left\langle x(t)-x^*+\frac{t}{\alpha-1}\dot{x}(t),-\frac{t}{\alpha-1}\nabla_x\tilde{f}(x(t),\mu(t))\right\rangle
\\= & \frac{4}{\alpha-1}t\left(\tilde{f}(x(t),\mu(t))-\tilde{f}(x^*,\mu(t))+2\kappa\mu(t)\right)
\\ & -2t\left\langle x(t)-x^*,\nabla_x\tilde{f}(x(t),\mu(t))\right\rangle.
\end{split}
\end{equation}
From the convexity of $\tilde{f}(x,\mu)$ with respect to $x$ for any fixed $\mu$, \eqref{fkf} and $\alpha\geq3$, we then deduce that
\begin{equation*}
\begin{split}
\frac{d}{dt}\mathcal{E}_{g,T}(t)\leq & \frac{2(3-\alpha)}{\alpha-1}t\left(\tilde{f}(x(t),\mu(t))-\tilde{f}(x^*,\mu(t))+2\kappa\mu(t)\right)+4\kappa t\mu(t)
\\ \leq & 4\kappa t\mu(t).
\end{split}
\end{equation*}
Integrating the above inequality from $t_0$ to $t$ yields
\begin{equation}\label{ui}
\mathcal{E}_{g,T}(t)-\mathcal{E}_{g,T}(t_0)\leq\int_{t_0}^{t}4\kappa\tau\mu(\tau)d\tau\leq\int_{t_0}^{\infty}4\kappa\tau\mu(\tau)d\tau<\infty,
\end{equation}
which uses condition $(H_1)$.
\\Recalling the definition of $\mathcal{E}_{g,T}$, we find
\begin{equation}\label{q22}
\begin{split}
& \frac{2}{\alpha-1}t^2\left(\tilde{f}(x(t),\mu(t))-\tilde{f}(x^*,\mu(t))+2\kappa\mu(t)\right)+(\alpha-1)\left\|x(t)-x^*
+\frac{t}{\alpha-1}\dot{x}(t)\right\|^2
 \\ \leq & C_2+2\int_{t_0}^{t}\tau\left\langle x(\tau)-x^*+\frac{\tau}{\alpha-1}\dot{x}(t),g(\tau)\right\rangle d\tau,
\end{split}
\end{equation}
where $C_2=\int_{t_0}^{\infty}4\kappa\tau\mu(\tau)d\tau+\frac{2}{\alpha-1}t_0^2\left(\tilde{f}(x(t_0),\mu(t_0))-\tilde{f}(x^*,\mu(t_0))+2\kappa\mu(t_0)\right)+(\alpha-1)\left\|x(t_0)-x^*
+\frac{t_0}{\alpha-1}\dot{x}(t_0)\right\|^2$ is a constant by \eqref{ui}. Then, from \eqref{fkf} and \eqref{q22}, we know
\begin{equation*}
\begin{split}
& \frac{1}{2}\left\|x(t)-x^*+\frac{t}{\alpha-1}\dot{x}(t)\right\|^2
\\ \leq & \frac{C_2}{2(\alpha-1)}+\frac{1}{\alpha-1}\int_{t_0}^{t}\left\|x(\tau)-x^*+\frac{\tau}{\alpha-1}\dot{x}(t)\right\|\left\|\tau g(\tau)\right\|d\tau.
\end{split}
\end{equation*}
Applying Lemma \ref{ap4}, we have
\begin{equation*}
\left\|x(t)-x^*+\frac{t}{\alpha-1}\dot{x}(t)\right\|\leq\left(\frac{C_2}{\alpha-1}\right)^{\frac{1}{2}}+\frac{1}{\alpha-1}\int_{t_0}^{t}\tau\|g(\tau)\|d\tau,
\end{equation*}
which together with $\int_{t_0}^{\infty}t\|g(t)\|dt<\infty$ gives
\begin{equation}\label{con}
\sup_{t\geq t_0}\left\|x(t)-x^*+\frac{t}{\alpha-1}\dot{x}(t)\right\|<\infty.
\end{equation}
Returning to \eqref{q22}, applying the Cauchy-Schwart inequality and by $\alpha>3$, we see that
\begin{equation*}
\begin{split}
& \frac{2}{\alpha-1}t^2\left(\tilde{f}(x(t),\mu(t))-\tilde{f}(x^*,\mu(t))+2\kappa\mu(t)\right)
\\ \leq & C_2+2\sup_{t\geq t_0}\left\|x(t)-x^*+\frac{t}{\alpha-1}\dot{x}(t)\right\|\int_{t_0}^{\infty}\|\tau g(\tau)\|d\tau<\infty.
\end{split}
\end{equation*}
In view of \eqref{fkf}, we deduce that
\begin{equation*}
f(x(t))-f(x^*)=O\left(\frac{1}{t^2}\right).
\end{equation*}
(ii) From \eqref{con}, we know
 \begin{equation*}
\begin{split}
\mathcal{E}_{g}(t):=& \frac{2}{\alpha-1}t^2\left(\tilde{f}(x(t),\mu(t))-\tilde{f}(x^*,\mu(t))+2\kappa\mu(t)\right)
\\ & +(\alpha-1)\left\|x(t)-x^*+\frac{t}{\alpha-1}\dot{x}(t)\right\|^2
\\ & +2\int_{t}^{\infty}\tau\left\langle x(\tau)-x^*+\frac{\tau}{\alpha-1}\dot{x}(t),g(\tau)\right\rangle d\tau
\end{split}
\end{equation*}
is well defined. Similar to the calculation to $\mathcal{E}_{g,T}$, we have that
\begin{equation}\label{qu}
\frac{d}{dt}\mathcal{E}_{g}(t)+2\frac{\alpha-3}{\alpha-1}t\left(\tilde{f}(x(t),\mu(t))-\tilde{f}(x^*,\mu(t))+2\kappa\mu(t)\right)\leq4\kappa t\mu(t).
\end{equation}
By integrating \eqref{qu} on $[t_0,t]$, we obtain
\begin{equation*}
\begin{split}
& \mathcal{E}_g(t)+2\frac{\alpha-3}{\alpha-1}\int_{t_0}^{t}\tau\left(\tilde{f}(x(\tau),\mu(\tau))-\tilde{f}(x^*,\mu(\tau))+2\kappa\mu(\tau)\right)d\tau
\\\leq & \int_{t_0}^{t}4\kappa\tau\mu(\tau)d\tau+\mathcal{E}_g(t_0).
\end{split}
\end{equation*}
Recalling the definition of $\mathcal{E}_g$ and by \eqref{fkf}, we obtain
\begin{equation*}
\begin{split}
& 2\int_{t}^{\infty}\tau\left\langle x(\tau)-x^*+\frac{\tau}{\alpha-1}\dot{x}(\tau),g(\tau)\right\rangle d\tau+2\frac{\alpha-3}{\alpha-1}\int_{t_0}^{t}\tau(\tilde{f}(x(\tau),\mu(\tau))
\\ & -\tilde{f}(x^*,\mu(\tau))+2\kappa\mu(\tau))d\tau
\leq\int_{t_0}^{t}4\kappa\tau\mu(\tau)d\tau+\mathcal{E}_g(t_0).
\end{split}
\end{equation*}
Rearranging the above inequality and using the Cauchy-Schwart inequality, we infer that
\begin{equation*}
\begin{split}
& 2\frac{\alpha-3}{\alpha-1}\int_{t_0}^{t}\tau\left(\tilde{f}(x(\tau),\mu(\tau))-\tilde{f}(x^*,\mu(\tau))+2\kappa\mu(\tau)\right)d\tau
\\\leq & \mathcal{E}_g(t_0)+2\sup_{t\geq t_0}\left\|x(t)-x^*+\frac{t}{\alpha-1}\dot{x}(t)\right\|\int_{t_0}^{\infty}\|\tau g(\tau)\|d\tau+\int_{t_0}^{t}4\kappa\tau\mu(\tau)d\tau.
\end{split}
\end{equation*}
Recalling $\alpha>3$ and \eqref{con}, under the condition $(H_1)$ and by $\int_{t_0}^{\infty}t\|g(t)\|dt<\infty$, we conclude that
\begin{equation}\label{xi}
\int_{t_0}^{\infty}\tau\left(\tilde{f}(x(\tau),\mu(\tau))-\tilde{f}(x^*,\mu(\tau))+2\kappa\mu(\tau)\right)d\tau<\infty.
\end{equation}
According to \eqref{fkf}, we obtain
\begin{equation*}
\int_{t_0}^{\infty}t\left(f(x(t)-{\rm min} f\right)dt<\infty,
\end{equation*}
which proves \eqref{q30}.
\\Next, let us show that \eqref{q20}. Taking the scalar product of \eqref{prob3} with $t^2\dot{x}(t)$, we obtain
\begin{equation*}
t^2\left\langle\ddot{x}(t),\dot{x}(t)\right\rangle+\alpha t\|\dot{x}(t)\|^2+t^2\left\langle\nabla_x\tilde{f}(x(t),\mu(t)),\dot{x}(t)\right\rangle=t^2\left\langle g(t),\dot{x}(t)\right\rangle.
\end{equation*}
Using the Chain rule, Cauchy-Schwart inequality and \eqref{fkg}, we get
\begin{equation}\label{asd}
\begin{split}
& \frac{1}{2}t^2\frac{d}{dt}\|\dot{x}(t)\|^2+\alpha t\|\dot{x}(t)\|^2+t^2\frac{d}{dt}\left(\tilde{f}(x(t),\mu(t))-\tilde{f}(x^*,\mu(t)) +2\kappa\mu(t)\right)
\\ \leq & \|tg(t)\|\|t\dot{x}(t)\|.
\end{split}
\end{equation}
\\Let us integrate \eqref{asd} on $[t_0,t]$, then we have
\begin{equation*}
\begin{split}
& \frac{t^2}{2}\|\dot{x}(t)\|^2-\frac{{t_0}^2}{2}\|\dot{x}(t_0)\|^2+(\alpha-1)\int_{t_0}^{t}s\|\dot{x}(s)\|^2ds+t^2(\tilde{f}(x(t),\mu(t))-\tilde{f}(x^*,\mu(t))
\\ & +2\kappa\mu(t))-t_0^2(\tilde{f}(x(t_0),\mu(t_0))-\tilde{f}(x^*,\mu(t_0))+2\kappa\mu(t_0))-2\int_{t_0}^{t}s(\tilde{f}(x(s),\mu(s))
\\ & -\tilde{f}(x^*,\mu(s))+2\kappa\mu(s))ds\leq\int_{t_0}^{t}\|sg(s)\|\|s\dot{x}(s)\|ds.
\end{split}
\end{equation*}
Using \eqref{fkf} again, we have
\begin{equation}\label{sim}
\begin{split}
\frac{1}{2}\|t\dot{x}(t)\|^2+(\alpha-1)\int_{t_0}^{t}s\|\dot{x}(s)\|^2ds\leq & C_3+2\int_{t_0}^{t}s(\tilde{f}(x(s),\mu(s))-\tilde{f}(x^*,\mu(s))
\\ & +2\kappa\mu(s))ds+\int_{t_0}^{t}\|sg(s)\|\|s\dot{x}(s)\|ds
\end{split}
\end{equation}
for some constant $C_3=\frac{{t_0}^2}{2}\|\dot{x}(t_0)\|^2+t_0^2(\tilde{f}(x(t_0),\mu(t_0))-\tilde{f}(x^*,\mu(t_0))+2\kappa\mu(t_0))$ only depending on the information at $t_0$. \\
From \eqref{xi} and $\alpha>1$, it follows that
\begin{equation*}
\frac{1}{2}\|t\dot{x}(t)\|^2\leq C_4+\int_{t_0}^{t}\|sg(s)\|\|s\dot{x}(s)\|ds,
\end{equation*}
where $C_4=C_3+\int_{t_0}^{\infty}\tau\left(\tilde{f}(x(\tau),\mu(\tau))-\tilde{f}(x^*,\mu(\tau))+2\kappa\mu(\tau)\right)d\tau$ is a constant.
Applying Lemma \ref{ap4}, we obtain
\begin{equation*}
\|t\dot{x}(t)\|\leq\left(2C_4\right)^{\frac{1}{2}}+\int_{t_0}^{t}\|sg(s)\|ds,
\end{equation*}
by $\int_{t_0}^{\infty}\|sg(s)\|ds<\infty$, which implies
\begin{equation}\label{xi1}
\sup _{t\geq t_0}\|t\dot{x}(t)\|<\infty.
\end{equation}
Returning to \eqref{sim}, it gives
\begin{equation*}
\begin{split}
(\alpha-1)\int_{t_0}^{t}s\|\dot{x}(s)\|^2ds\leq& C_3+2\int_{t_0}^{t}s\left(\tilde{f}(x(s),\mu(s))-\tilde{f}(x^*,\mu(s))+2\kappa\mu(s)\right)ds
\\
&+\sup _{t\geq t_0}\|t\dot{x}(t)\|\int_{t_0}^{\infty}\|sg(s)\|ds,
\end{split}
\end{equation*}
which implies \eqref{q20}.
\\Combining \eqref{con} with \eqref{xi1}, we get
\begin{equation}\label{xi2}
\sup_{t\geq t_0}\|x(t)\|<\infty.
\end{equation}
Recalling the definition of $E$ in \eqref{num}, by \eqref{q20} and \eqref{xi}, we deduce that
\begin{equation*}
\int_{t_0}^{\infty}tE(t)dt<\infty.
\end{equation*}
Same as the proof of (ii) in Theorem \ref{th} we can deduce that
\begin{equation*}
f(x(t))-\min f=o\left(\frac{1}{t^2}\right).
\end{equation*}
\qed
\subsection{Weak convergence of trajectory under perturbations}
\label{sec:3}
\begin{theorem}\label{}
 Assume ${\rm argmin} f\neq\emptyset$, $\int_{t_0}^{\infty}t\|g(t)\|dt<\infty$ and $(H_1)$ is ture. Let $x: [t_0,\infty)\rightarrow\mathcal{H}$ be the trajectory of \eqref{prob3} with $\alpha>3$. Then $x(t)$ converges weakly to an element of {\rm argmin} $f$, as $t\rightarrow\infty$.
\end{theorem}
{\it Proof}\quad
Applying Cauchy-Schwarz inequality in \eqref{q17}, we find that
\begin{equation*}
\begin{split}
& t\ddot{h}(t)+\alpha\dot{h}(t)+t\left(\tilde{f}(x(t),\mu(t))-\tilde{f}(z,\mu(t))+2\kappa\mu(t)\right)
\\ \leq & t\|\dot{x}(t)\|^2+\|x(t)-z\|\|tg(t)\|+2\kappa t\mu(t),
\end{split}
\end{equation*}
where $h$ is defined in \eqref{num1}.
\\By virtue of \eqref{fkf} and letting $z=x^*\in\argmin f$, we have
\begin{equation*}
\begin{split}
t\ddot{h}(t)+\alpha\dot{h}(t)\leq l(t),
\end{split}
\end{equation*}
where $l(t):=t\|\dot{x}(t)\|^2+\left(\sup_{t\geq t_0}\|x(t)-x^*\|\right)\|tg(t)\|+2\kappa t\mu(t)\geq0$. Under $(H_1)$, \eqref{q20}, \eqref{xi2} and $\int_{t_0}^{\infty}t\|g(t)\|dt<\infty$, we deduce that $l(t)\in L^1(t_0,\infty)$. By Lemma \ref{ap3}, we know $\lim_{t\rightarrow\infty}h(t)$ exists, where $x^*$ can be any element in argmin$f$. Moreover, let $\bar{x}$ be a sequential limit point of $x(t)$ on $[t_0,\infty)$ with convergence sequence $\{t_n\}$. Similar to the proof of Theorem \ref{theo}, using Proposition \ref{soul1} we obtain that $\bar{x}\in {\rm argmin}f$. Hence, the proof is completed by Lemma \ref{ap1}.
\qed
\section{Numerical experiments}\label{s5}
\begin{figure}[H]
\centerline{
 \subfigure[]{\includegraphics[width=0.53 \textwidth]{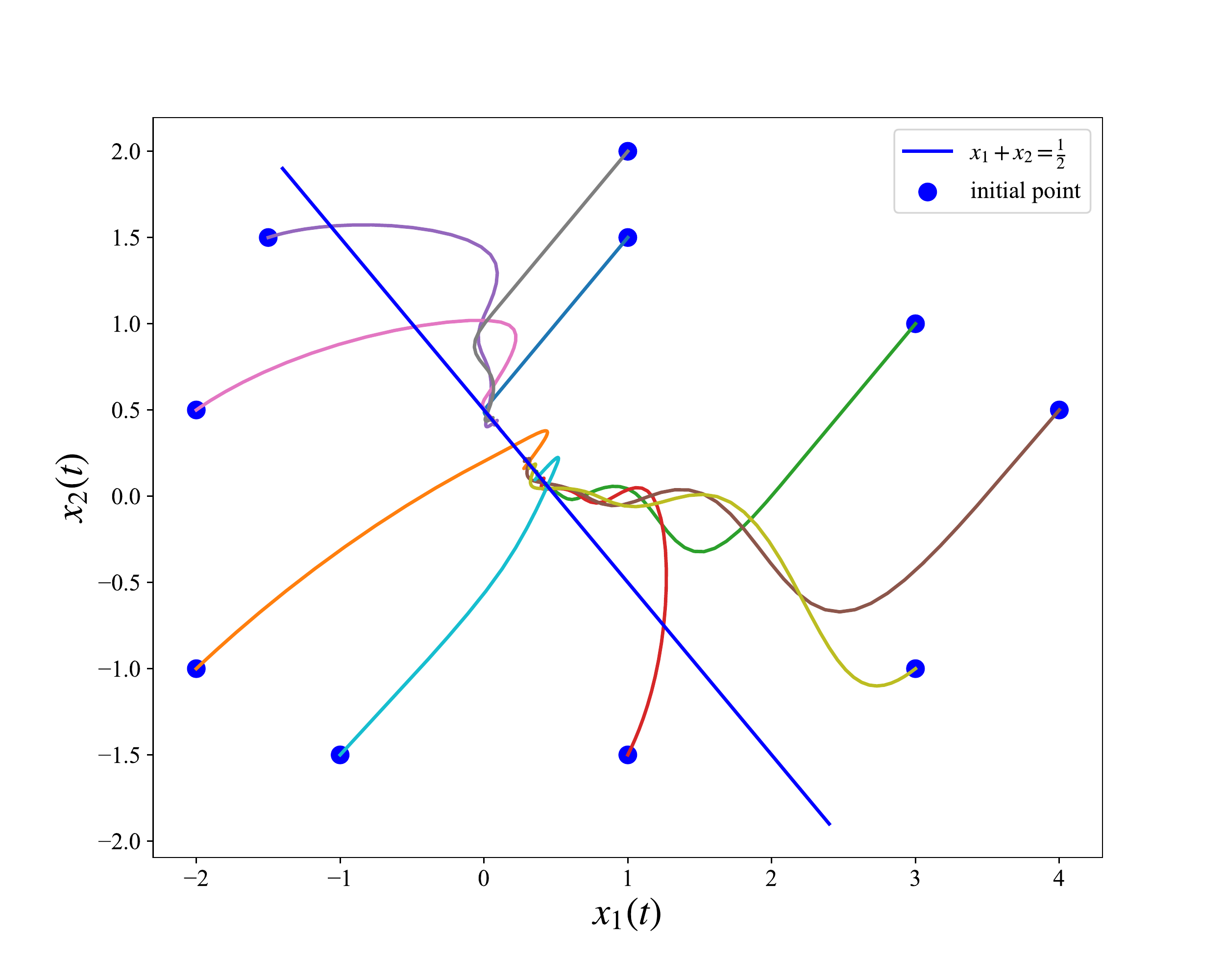}}
    \subfigure[]{\includegraphics[width=0.53 \textwidth]{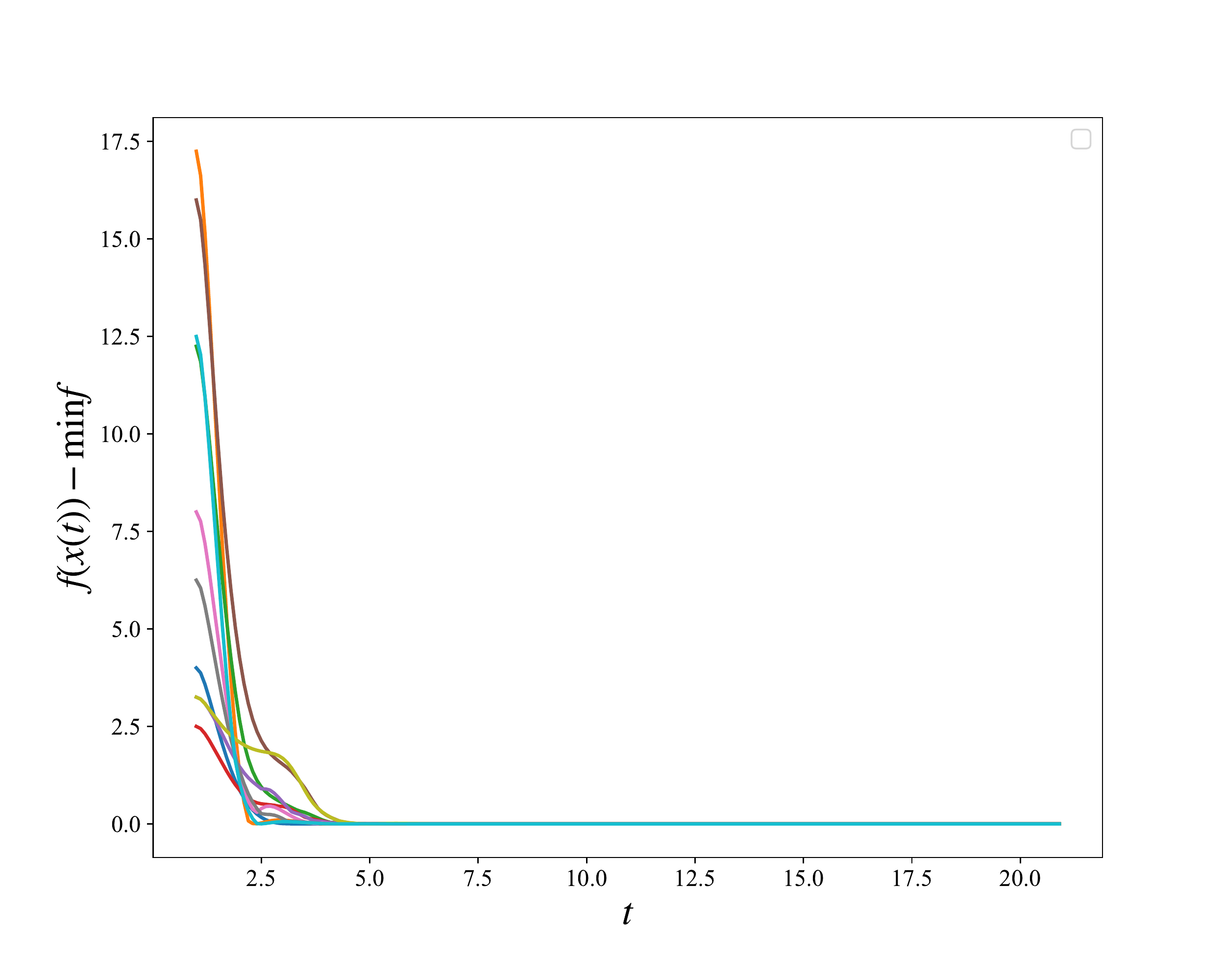}}
}
  \caption{Convergence of trajectories and function values associated with dynamic algorithm \eqref{prob2} for Example \ref{exa2}}\label{pic1}
\end{figure}

\begin{figure}[H]
\centerline{
 \subfigure[]{\includegraphics[width=0.53 \textwidth]{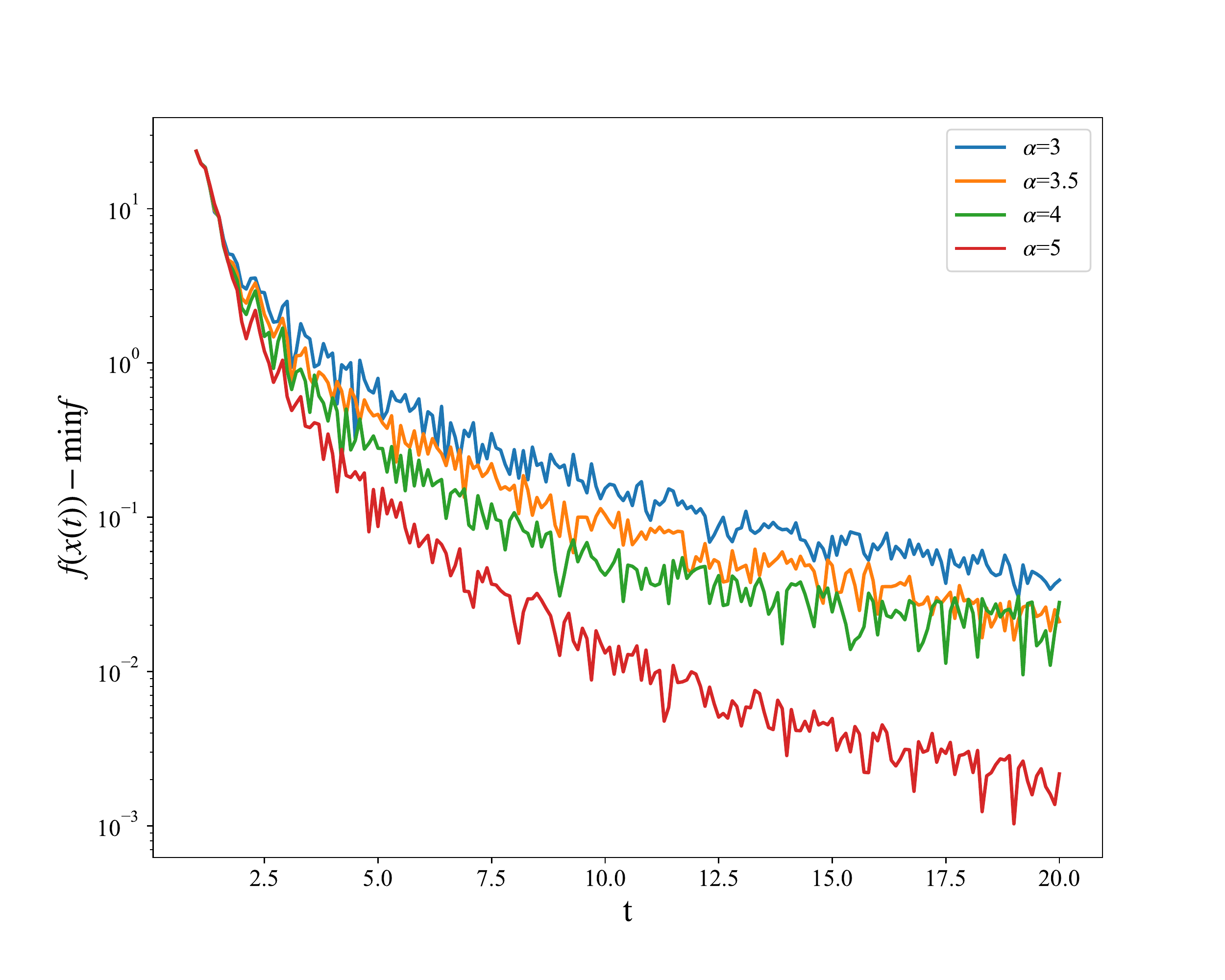}}
    \subfigure[]{\includegraphics[width=0.53 \textwidth]{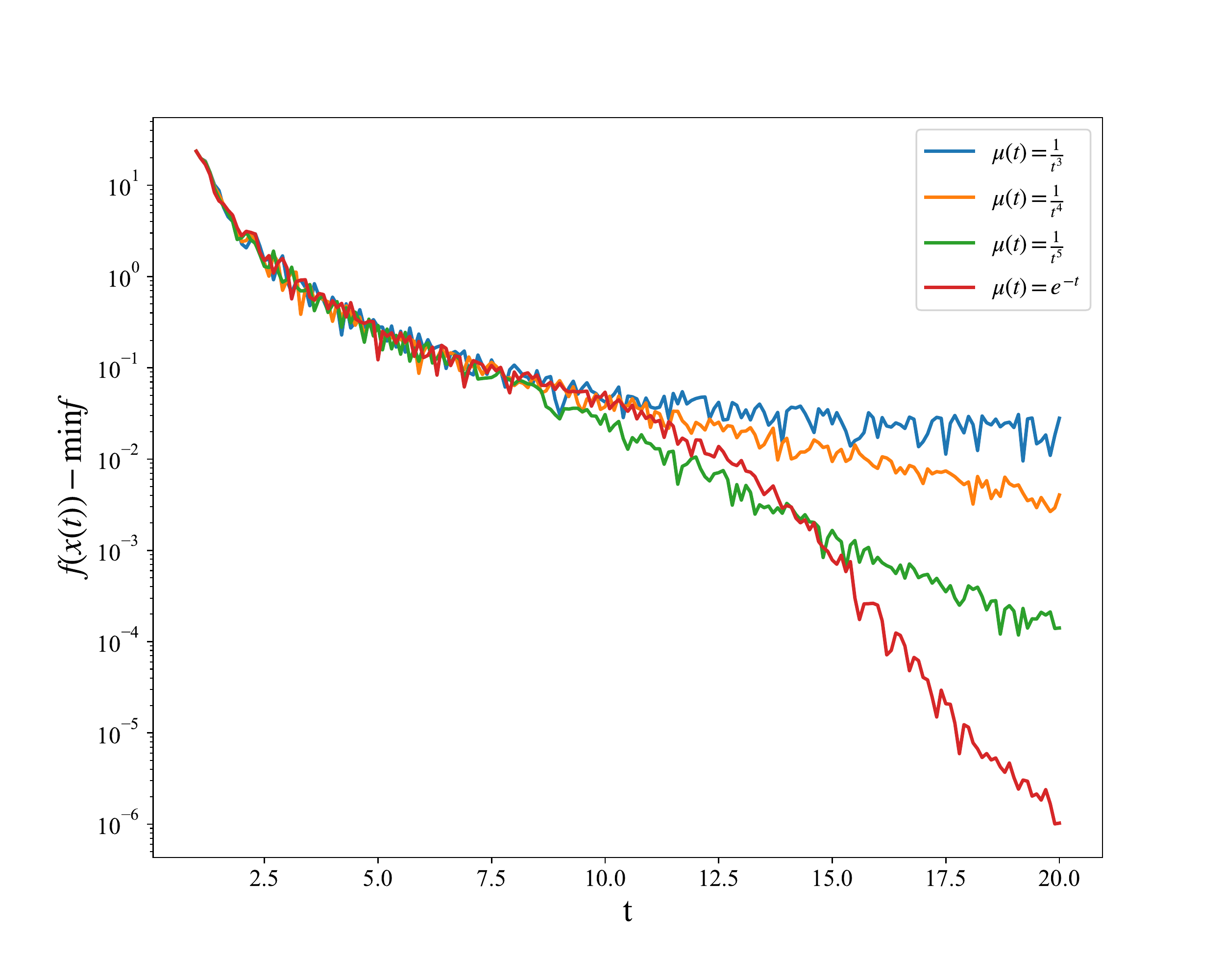}}
}
  \caption{Convergence of function values associated with dynamic algorithm \eqref{prob2} for Example \ref{exa8}}\label{pic8}
\end{figure}

\begin{figure}[H]
\centerline{
 {\includegraphics[width=0.53 \textwidth]{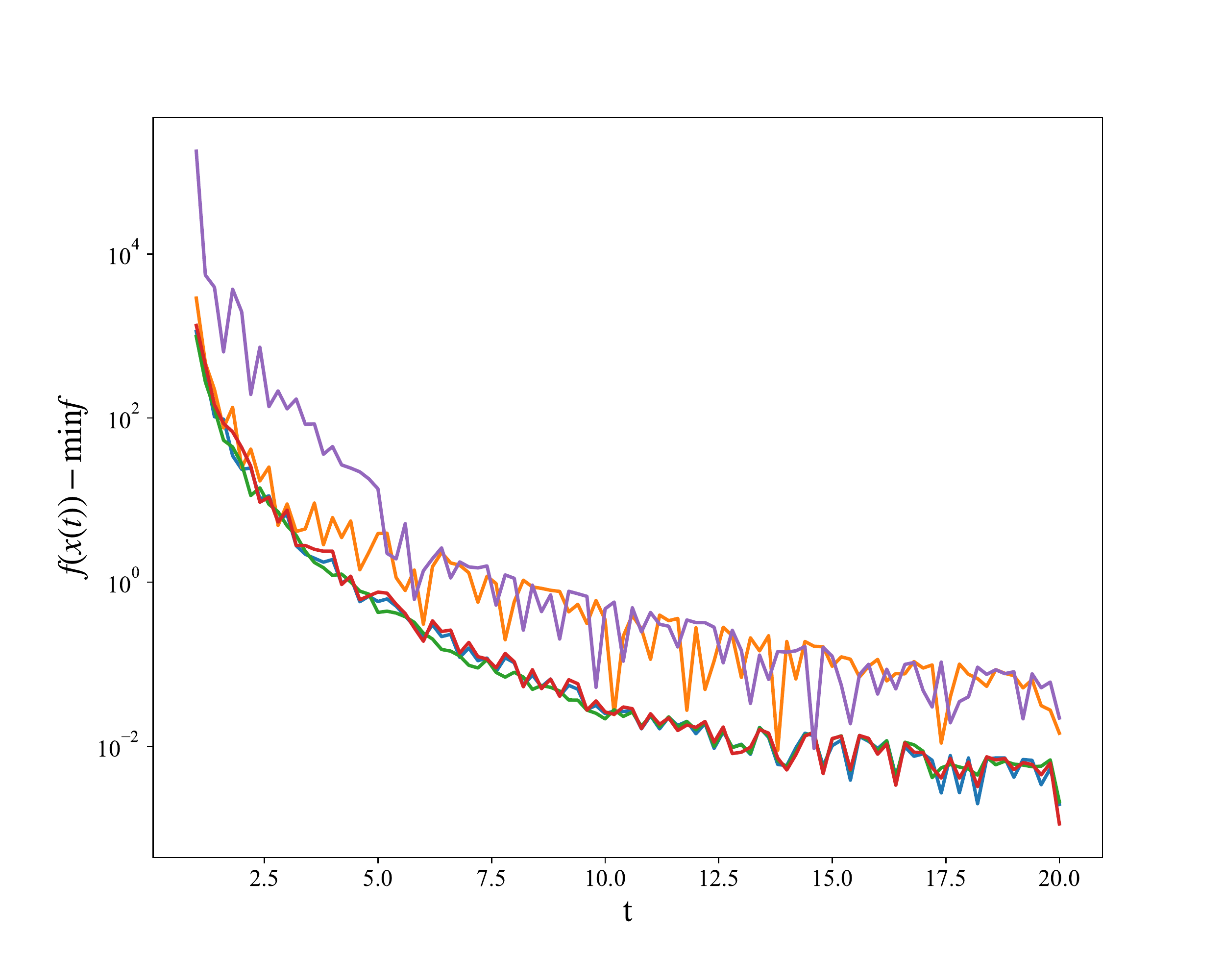}}
}
  \caption{Convergence of function values associated with dynamic algorithm \eqref{prob2} for Example \ref{exa3}}\label{pic3}
\end{figure}

In this section, we report three numerical experiments to verify the theoretical results of dynamic algorithms \eqref{prob2} and \eqref{prob3}. All experiments are performed in Python 3.7.3 on a Lenovo PC (2.30GHz, 8.00GB of RAM).
\begin{example}\label{exa2}
Consider the following nonsmooth convex optimization problem in $\mathbb{R}^2$, namely,
\begin{equation}\label{ex2}
\min f(x_1, x_2)=(x_1+x_2-1)^2+|x_1|+\max \{x_2,0\}.
\end{equation}
\end{example}
 We can deduce that optimal solution set of \eqref{ex2} is $\hat{X}:=\{x:x_1+x_2=\frac{1}{2},x_1\geq0, x_2\geq0\}$, and the optimal value is $f(\hat{x})=\frac{3}{4},\forall \hat{x}\in \hat{X}$.
  \\ Let $\alpha=7$ and $\mu(t)=\frac{1}{t^3}$. Fig. \ref{pic1} illustrates that trajectories of \eqref{prob2} with ten random initial points converge to some elements in $\hat{X}$, and shows the convergence rate on objective values.
\begin{figure}[H]
\centerline{
 \subfigure[]{\includegraphics[width=0.53 \textwidth]{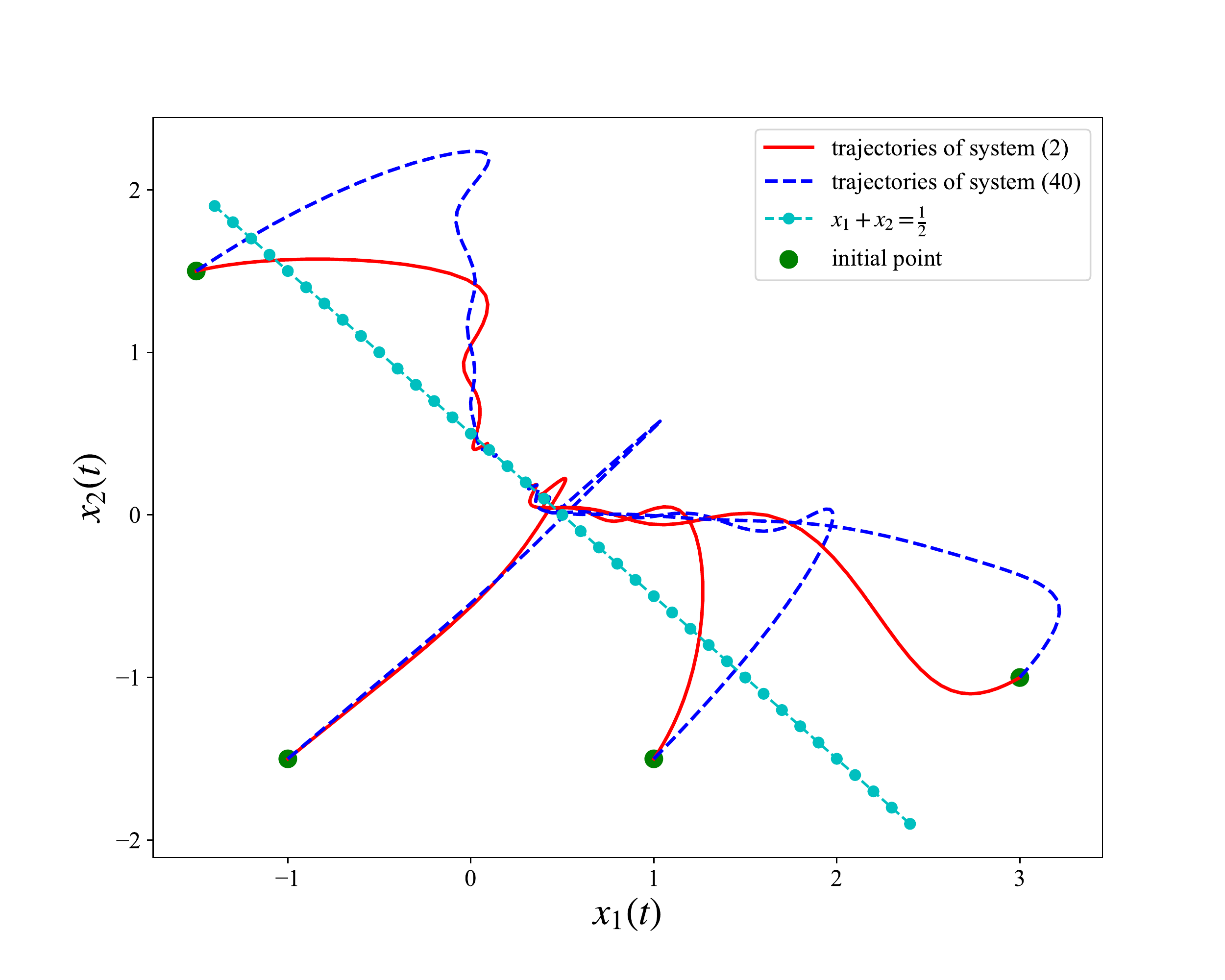}}
    \subfigure[]{\includegraphics[width=0.53 \textwidth]{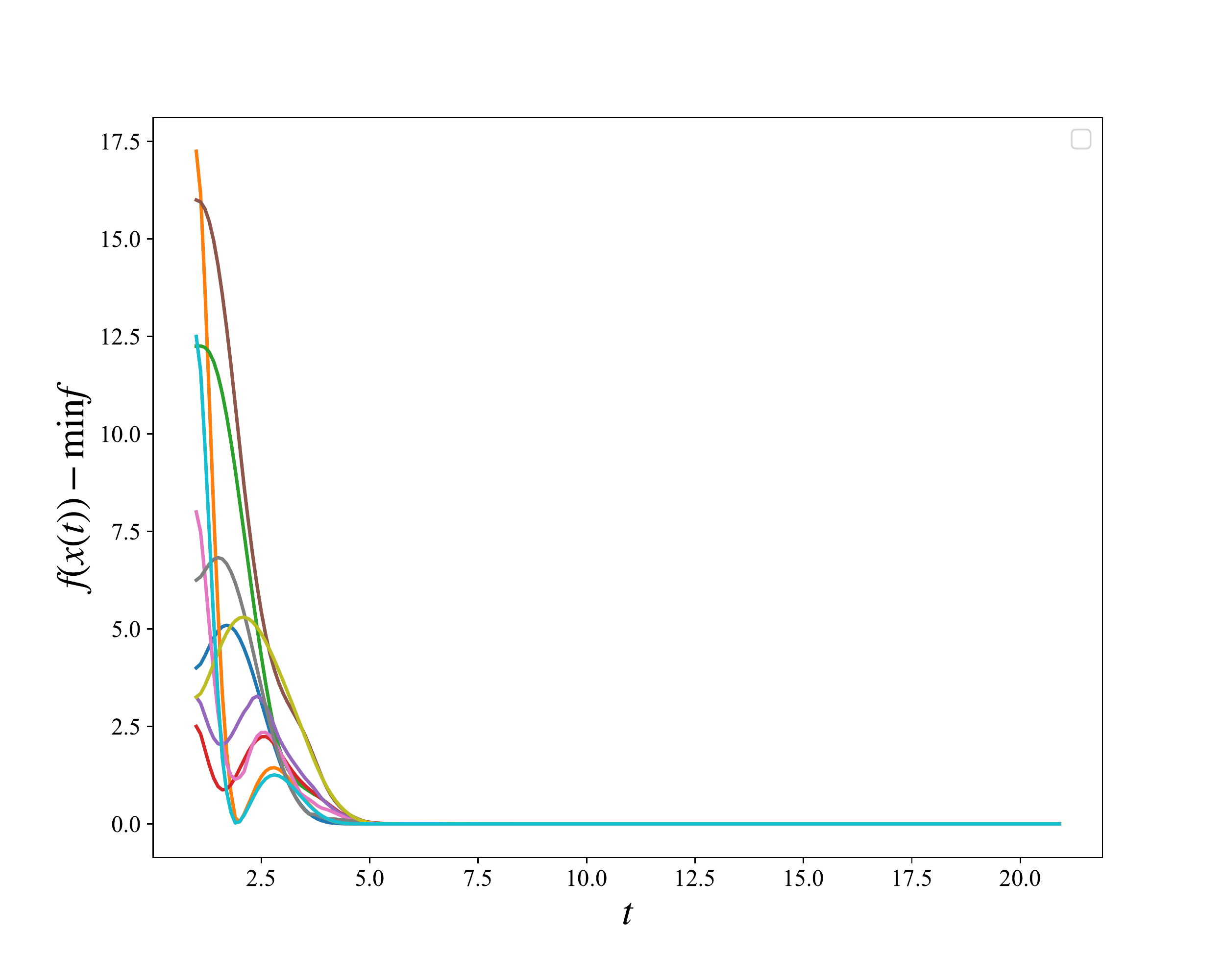}}
}
  \caption{Convergence of trajectories and function values associated with dynamic algorithm \eqref{prob3} for Example \ref{exa2}}\label{pic2}
\end{figure}
For the above Example \ref{exa2}, we take as a perturbation $g(t)=20e^{-t}$ to verify the theoretical results of dynamic algorithm \eqref{prob3}. The corresponding results are presented in Fig. \ref{pic2}-(a) with four random initial points, from which we can see the stability of dynamic algorithm \eqref{prob2} under the perturbation. And Fig. \ref{pic2}-(b) shows the convergence rate on objective values of dynamic algorithm \eqref{prob3} with ten random initial points.
\begin{example}\label{exa8}
We consider the following nonsmooth convex optimization problem,
\begin{equation}
\min f(x):=\|Ax-b\|^2_2+\|Dx-d\|_1,
\end{equation}
where $A\in \mathbb{R}^{20\times10}$, $b\in \mathbb{R}^{20}$, $D\in \mathbb{R}^{50\times10}$ and $d\in \mathbb{R}^{50}$ are randomly generated as follows:
$$A=\rm{randn}(20,10);D=\rm{randn}(50,10);$$
$$x^*=\rm{randn}(10,1);b=A*x^*;d=D*x^*.$$
\end{example}
From the above generation, we know that $x^*$ is an optimal solution of Example \ref{exa8} and its optimal value is 0.
\\$\bullet$ Let $\mu(t)=\frac{1}{t^3}$. Fig. \ref{pic8}-(a) shows the influence of $\alpha$ on the convergence rate of function values of dynamic algorithm \eqref{prob2} with the same initial value. It can be seen from Fig. \ref{pic8}-(a) that for different values of $\alpha$, each objective function value converges to the optimal value along the trajectory of dynamic algorithm \eqref{prob2}, and the larger of $\alpha$, the faster the convergence rate of the function value.
\\$\bullet$ Let $\alpha=4$. Fig. \ref{pic8}-(b) shows the influence of $\mu(t)$ on the convergence rate of function values of dynamic algorithm \eqref{prob2} with the same initial value. We can see that for selecting different $\mu(t)$, each objective function value also converges to the optimal value along the trajectory of dynamic algorithm \eqref{prob2}, and the faster $\mu(t)$ decreases as $t\rightarrow\infty$, the faster the convergence rate of the function value.
\begin{example}\label{exa3}
For example \ref{exa8}, we consider the high-dimensional case
\begin{equation}
\min f(x):=\|Ax-b\|^2_2+\|Dx-d\|_1,
\end{equation}
where $A\in \mathbb{R}^{200\times100}$, $b\in \mathbb{R}^{200}$, $D\in \mathbb{R}^{500\times100}$ and $d\in \mathbb{R}^{500}$ are randomly generated as follows:
$$A=\rm{randn}(200,100);D=\rm{randn}(500,100);$$
$$x^*=\rm{randn}(100,1);b=A*x^*;d=D*x^*.$$
\end{example}
 Fig. \ref{pic3} presents the fast convergence rate on the function values of dynamic algorithm \eqref{prob2} with five random initial points.
\section{Conclusions}
In this paper, we focused on the asymptotic convergence of dynamic algorithm \eqref{prob2} and its a perturbed version \eqref{prob3} for solving convex optimization problem \eqref{x1}, where the smoothing method is used to overcome the gradient Lipschitz condition of the objective function. Firstly, we used Cauchy-Lipschitz-Picard theorem to prove the global existence and uniqueness of the trajectory of dynamic algorithm \eqref{prob2}. Then by constructing an appropriate energy functions, we showed that the convergence rate on the objective values is $O\left(t^{-2}\right)$ as $\alpha\geq3$, and $o\left(t^{-2}\right)$ as $\alpha>3$, which are same as the results of dynamic algorithm \eqref{in1} for solving the corresponding continuous differentiable convex optimization problems. In addition, we proved that the trajectory of \eqref{prob2} is weakly convergent to an optimal solution of problem \eqref{x1}. For the perturbed second-order dynamic algorithm \eqref{prob3}, we verified that it has the same convergence properties as \eqref{prob2} under a proper condition on the perturbation. Finally, we illustrated the theoretical results by some numerical examples.

\vspace{0.5cm}

\textbf{Funding}
This work is funded by the National Science Foundation of China (No: 11871178).

\vspace{0.5cm}

\textbf{Data and code availability} The data and code that support the fndings of this study are available from the corresponding author upon request.

\end{document}